\documentclass[12pt]{amsart}
\usepackage{amssymb,amsmath,graphics,verbatim}
\usepackage{latexsym}
\usepackage{eucal}
\usepackage{a4wide}
\usepackage{color}

\newtheorem{theorem}{Theorem}[section]
\newtheorem{lemma}[theorem]{Lemma}
\newtheorem{remark}[theorem]{Remark}
\newtheorem{prop}[theorem]{Proposition}
\newtheorem{cor}[theorem]{Corollary}

\def \N {\mathbb{N}}
\def \R{\mathbb{R}}
\def \C{\mathbb{C}}

\def \supp{\rm supp}

\def \ph{\varphi}
\def \Om{\Omega}

\def \grad{\nabla}

\def \S{\mathcal{S}}

\def \t{\tilde}
\def \half{\frac{1}{2}}

\numberwithin{equation}{section}

\begin{document}
\title{The $L^p$ Carleman Estimate and a Partial Data Inverse Problem}

\author[Chung]{Francis J. Chung}
\address{Department of Mathematics, University of Kentucky, Lexington, KY, USA}
\email{fj.chung@uky.edu}

\author[Tzou]{Leo Tzou}
\address{School of Mathematics and Statistics, University of Sydney, Sydney, Australia}
\email{leo@maths.usyd.edu.au}



\begin{abstract}
We construct an explicit Green's function for the conjugated Laplacian $e^{-\omega \cdot x/h}\Delta e^{-\omega \cdot x/h}$, which let us control our solutions on roughly half of the boundary. We apply the Green's function to solve a partial data inverse problem for the Schr\"odinger equation with potential $q \in L^{n/2}$.  We also use this Green's function to derive $L^p$ Carleman estimates similar to the ones in Kenig-Ruiz-Sogge \cite{krs}, but for functions with support up to part of the boundary. 
\end{abstract}

\maketitle

\begin{section}{Introduction}
In this article we give an explicit construction of a ``Dirichlet Green's function" for the conjugated Laplacian $e^{-x\cdot\omega/h} h^2\Delta e^{x\cdot\omega/h}$ on a bounded smooth domain $\Omega \subset \R^n$ for $n\geq 3$. This Green's function immediately gives various ($L^2$ and $L^p$) Carleman estimates similar to those in Kenig-Ruiz-Sogge \cite{krs} and Kenig-Sj\"ostrand-Uhlmann \cite{ksu} (linear weight case), %\fc{Up to change of vars, log weight should be linear too?} 
for functions in $C^\infty(\Omega)$ with nontrivial boundary conditions. We also apply the Green's function to solve the partial data inverse Schr\"{o}dinger problem with unbounded potential in $L^{n/2}(\Omega)$ for $n\geq 3$.

The main result is the construction of the Green's function. Let $\omega\in \R^n$ be a unit vector and let $\Gamma\subset\partial\Omega$ be an open subset which is compactly contained in $\{x\in\partial\Omega \mid\nu(x)\cdot \omega >0\}$. If $p' = \frac{2n}{n+2} <2< p = \frac{2n}{n-2}$, we have the following theorem, proved by an explicit construction via heat flow.
\begin{theorem}
\label{green's function}
Suppose $h > 0$ is sufficently small. Then there exists an operator $G_\Gamma : L^{p'}(\Omega) \to L^p(\Omega)$ which satisfies \[ e^{-x\cdot\omega/h} h^2\Delta e^{x\cdot\omega/h} G_\Gamma = I\] and the estimates
\[\|G_\Gamma\|_{L^2\to H^1} \leq Ch^{-1},\ \ \|G_\Gamma\|_{L^{p'}\to L^p} \leq Ch^{-2}.\]
Furthermore, for all $f\in L^{p'}$, $G_\Gamma f \in H^1(\Omega)$ and $G_\Gamma f \mid_{\Gamma} = 0$.
\end{theorem}
We use the Green's function to prove the following Carleman estimates. Let $H^1(\Omega)$ denote the semiclassical Sobolev space. Define $H^1_\Gamma(\Omega)\subset H^1(\Omega)$ to be the space of functions with vanishing trace along $\Gamma$ and let $H^{-1}_\Gamma(\Omega)$ be its dual.  
\begin{theorem}
\label{carleman estimates}
Let $u \in C^2(\bar\Omega)$ be a function which vanishes along $\partial\Omega$ and $\partial_\nu u \mid_{\Gamma^c} = 0$. One then has the Carleman estimates
\[
\|u\|_{L^2(\Omega)} \leq \frac{C}{h}\|h^2 \Delta_\phi^{*} u\|_{H^{-1}_\Gamma(\Omega)} 
,\ 
 \|u\|_{L^p(\Omega)} \leq C\|\Delta_\phi^{*} u\|_{L^{p'}(\Omega)}
\]
for all $h>0$ sufficiently small.
\end{theorem}
\begin{remark} A modification of the argument presented here can also yield a boundary term of $h^{-1/p}\|\partial_\nu u \|_{L^{p'}(\Gamma)}$ on the left-side of the $L^p$ inequality. 
\end{remark}
%Note that in the ``shifted" $L^2$-Sobolev estimate the domain does not need to be enlarged as was necessary in \cite{ksu}. 
%\fc{I don't think KSU actually shift the Sobolev estimate at all, though DKSU magnetic Schrodinger paper does.} The first estimate is similar to the one in \cite{Chuthesis}, but the proof is much improved. 

The second estimate differs from other $L^p$ Carleman estimates for the Laplacian in that it allows for $u$ with nontrivial boundary conditions. 

Another application of this Green's function is the resolution of the partial data Calder\'on problem with unbounded potentials. Let $\Omega$ be a smooth domain contained in $\R^n$, with $n \geq 3$, and let $\omega_0\in \R^n$ be a unit vector. Define \[\Gamma_{\pm}^0 := \{x\in \partial\Omega \mid \pm\nu(x) \cdot \omega_0 \geq 0\}\]
and let ${\bf F} \subset \partial\Omega$ be an open neighbourhood containing $\Gamma_+^0$ and ${\bf B}\subset \partial\Omega$ be an open neighbourhood containing $\Gamma_-^0$.  

If zero is not an eigenvalue of the operator $-\Delta + q$, then $q \in L^{n/2}(\Omega)$ gives rise to a well-defined Dirichlet-to-Neumann map 
\[
\Lambda_q: H^{\half}(\partial \Om) \rightarrow H^{-\half}(\partial \Om). 
\]
(We refer the reader to the appendix of \cite{DosKenSal} for the definition of the Dirichlet-to-Neumann map for $q\in L^{n/2}(\Omega)$.) We have the following theorem. 

\begin{theorem}
\label{main theorem}
Let $q_1, q_2 \in L^{n/2}(\Omega)$ be such that $\Lambda_{q_1}f \mid_{{\bf F}} = \Lambda_{q_2}f \mid_{{\bf F}}$ for all $f\in C^\infty_0({\bf B})$. Then $q_1 = q_2$.
\end{theorem}

The regularity assumption that $q_j \in L^{n/2}$ is considered optimal in the context of well-posedness theory for the Dirichlet problem; it is also the optimal assumption for the strong unique continuation principle to hold (see \cite{JerKen} for more). 

We will provide some brief historical context for these theorems. The construction of the Green's function for the conjugated Laplace operator was established by Sylvester-Uhlmann \cite{su} using Fourier multipliers with characteristic sets. The authors proved an $L^2$ estimate for their Green's function and used it to solve the Calder\'{o}n problem in dimensions $n\geq 3$ for bounded potentials. Chanillo in \cite{ch} showed that the Sylvester-Uhlmann Green's function also satisfies an $L^p\to L^{p'}$ estimate by applying using the result of Kenig-Ruiz-Sogge \cite{krs}. This allowed Chanillo to solve the inverse Schr\"odinger problem with full data for small potentials in the Fefferman-Phong class (which contains $L^{n/2}$). Related results were also proved by Lavine-Nachman \cite{LavNac} and Dos Santos Ferreira-Kenig-Salo \cite{DosKenSal}. 

The drawback to the Fourier multiplier construction of the Green's function is that boundary conditions cannot be imposed. Bukhgeim-Uhlmann \cite{BukUhl} and Kenig-Sj\"ostrand-Uhlmann \cite{ksu} found a way to use Carleman estimates overcome this problem and prove results for the Calder\'{o}n problem with \emph{partial} boundary data. Due to its versatility and robustness, this technique has since become the standard tool for solving partial data elliptic inverse problems. The review article \cite{KenSalreview} contains an excellent overview of recent work in partial data Calder\'{o}n-type problems; examples for other elliptic inverse problems can be found in \cite{SalTzo}, \cite{SalTzo2}, \cite{KruLasUhl}, \cite{ChuSalTzo}, and \cite{ChuOlaSalTzo}.

The Carleman estimates in these papers are typically proved via an integration-by-parts procedure so that boundary conditions can be kept in check. The limitation of this approach is that only $L^2$-type estimates can be derived; none of the available techniques adapt well to $L^p$ setting for functions with boundary conditions. Thus for $q \notin L^{\infty}$, there are very few partial data results for the Calder\'{o}n problem -- see \cite{KruUhl32} for an example of what can be obtained by previous methods.

The Carleman estimate approach has the additional drawback that the Green's function one ``constructs" is an abstract object arising from general statements in functional analysis, like the Hahn-Banach or Riesz representation theorems. This makes partial data reconstruction procedures like the ones in \cite{NacStr} much more difficult to implement in a concrete setting than equivalent ones like \cite{Nac} for full data. %In a forthcoming article the authors will apply the Green's function constructed here to address this issue.

The Green's function we construct in Theorem \ref{green's function} has the explicit representation of the Fourier multiplier Green's function of Sylvester-Uhlmann while at the same time allowing the boundary control of the Carleman estimate approaches. Furthermore, due to its explicit representation as a parametrix, one can easily deduce $L^p$-type estimates as well as $L^2$-type estimates. 
%
%Another advantage of our Green's function is that it is essentially obtained by composing the Green's function of Sylvester-Uhlmann with a heat flow which one can explicitly solve. 
%
In a forthcoming article the authors intend to apply the Green's function constructed here to the problem of reconstruction. One expects that in the context of computational algorithms this Green's function would open the door to direct inversion methods for partial data Calder\'on problems in $n\geq 3$ which is parallel to the full data case examined in \cite{num2, num4, num1, num3}.

We give a brief exposition of our approach. The key observation is that there is a global $\Psi$DO factorization of the conjugated Laplacian $h^2 \Delta_\phi := e^{-\omega\cdot x/h}h^2 \Delta e^{\omega\cdot x/h}$ into an elliptic operator $J$ resembling a heat flow and a first-order operator $Q$ which has the same characteristic set as $h^2\Delta_\phi$. One can then construct an inverse for $J$ (and thus $h^2 \Delta_\phi$) with Dirichlet boundary conditions by solving the heat flow with zero initial condition.

This way of factoring $h^2\Delta_\phi$ is in the spirit of \cite{Chuthesis}. However, in our case the factorization is global and occurs on the level of symbols so there will be error terms and they pose a challenge in the construction of the parametrix. As such this necessitates a modified factorization which differs from that of \cite{Chuthesis} (see \eqref{operator factorization} and the discussions which follow) to obtain the suitable estimates for the remainders of the parametrix. 

This article is organized in the following way. In Section \ref{psido} we develop a $\Psi$DO calculus which is compatible with our symbol class - proofs are given in the appendix. In Section 3 we invert a heat flow in the context of this $\Psi$DO calculus and solve the Dirichlet problem for this heat flow. In section 4 we restate some facts about the Sylvester-Uhlmann Green's function in the semiclassical setting and derive a factorization for the operator $h^2\Delta_\phi$ involving the heat operator described in the previous section. In section 5 we use this factorization to construct a parametrix with Dirichlet boundary conditions, and in section 6 we turn the parametrix into a Dirichlet Green's function $G_\Gamma$ and prove Theorem \ref{carleman estimates}. Section 7 is devoted to proving Theorem \ref{main theorem} using complex geometric optics solutions constructed with the help of $G_\Gamma$.

\vspace{4mm}

\noindent \textbf{Acknowledgements}: The authors would like to thank the organizers of the Program on Inverse Problems at the Institut Henri Poincar\'{e}, where this project began. We would also like to thank Henrik Shahgholian of KTH and Yishao Zhou of Stockholm University for their hospitality during the summer of 2016. In addition, we would like to thank Boaz Haberman for several helpful discussions, and Sagun Chanillo for helping to explain the proof of Lemma \ref{krs estimate}.

\end{section}

\begin{section}{Elementary Semiclassical $\Psi$DO theory}\label{psido}
We collect a set of facts about semiclassical pseudodifferential operators and also use this opportunity to establish some notations and conventions which we will use throughout. Proofs are contained in the Appendix.  
\begin{subsection}{Mixed Sobolev Spaces}

In this article we define the semiclassical Sobolev spaces with the norm
\begin{eqnarray}
\label{sobolev}
\|u\|_{W^{k,r}_{scl}(\R^n)} := \| \langle hD\rangle^k  u\|_{L^r}.
\end{eqnarray}
For $k\in \N$ it turns out that this definition is equivalent to the one involving derivatives:
\[ \|u\|^r_{W^{k,r}_{scl}} = \sum\limits_{|\alpha |\leq k} \| (hD)^\alpha u\|^r_{L^r}.\]
(Hereafter we will drop the ``scl'' subscript: unless otherwise stated, all of our Sobolev spaces will be semiclassical.) 
Choose coordinates $(x',x_n)$ on $\R^n$, with $x' \in \R^{n-1}$ and $x_n \in \R$, and let $(\xi',\xi_n)$ be the corresponding coordinates on the cotangent space.
%The equivalence is trivial for the case when $k$ is even. For $k=1$ this can be easily seen by 
%{\Small\[\| \langle hD\rangle u\|_{L^r} = \| \langle hD\rangle^2 \langle hD\rangle^{-1} u \|_{L^r} \leq C\sum_{\substack{|\beta|\leq 1\\ |\alpha| = 1}} \| (hD)^\alpha  \langle hD\rangle^{-1}(hD)^\beta u\|_{L^r} \leq C\sum_{|\beta|\leq 1} \| (hD)^\beta u \|_{L^r}\]}
%The last inequality comes from the fact that $\xi_j \langle \xi\rangle^{-1}$ is a semiclassical Fourier multiplier on $L^r$ due to Mihlin's theorem. The equivalence for all odd $k$ then follows.
An immediate consequence of the norm equivalence stated above is that $\langle \xi' \rangle$ is a multiplier from $W^{1, r} (\R^n)\to L^r (\R^n)$. Indeed,
\begin{eqnarray}
\label{<hD'> map}
\| \langle hD'\rangle u \|^r_{L^r(\R^n)} &=& \int_{-\infty}^\infty  \| \langle hD'\rangle u(x',x_n)\|^r_{L^r_{x'}}dx_n \\
&\lesssim& \int_{-\infty}^\infty  \sum_{|\alpha|\leq 1}\| (hD')^\alpha u(x',x_n)\|^r_{L^r_{x'}}dx_n \leq  \sum_{|\alpha| = 1}\|(hD)^{\alpha} u\|^r_{L^r}. \nonumber 
\end{eqnarray}

Now define the mixed Sobolev norms for $u \in C^{\infty}_0$ by
\begin{eqnarray}
\label{mixed sobolev norms}
\|u\|_{W^{k,r}(\R^{n-1})W^{\ell, r}(\R^n)} := \| \langle hD'\rangle^k \langle hD\rangle^\ell u\|_{L^r}
\end{eqnarray}
and use these to define the mixed norm spaces $W^{k,r}(\R^{n-1})W^{\ell, r}(\R^n)$. For convenience we will drop the $\R^{n-1}$ and $\R^n$ in this notation and use the convention that the first $W^{k,r}$ denotes multiplication by $\langle hD'\rangle^k$ and the second $W^{\ell,r}$ denotes multiplication by $\langle hD\rangle^\ell$.

With this definition we have that for $k \geq 0$, \begin{eqnarray}\label{mixed space embedding}W^{-k,r} W^{\ell,r } \subset W^{l-k, r}(\R^n).\end{eqnarray} Indeed, one can write
\[ u = \langle hD'\rangle^k \langle hD\rangle^{-k} \langle hD\rangle^{-\ell + k} \langle hD'\rangle^{-k} \langle hD\rangle^\ell u\]
and use the fact that $\langle hD'\rangle^k \langle hD\rangle^{-k}$ is a multiplier on $L^r$ by \eqref{<hD'> map} and that \[\langle hD'\rangle^{-k} \langle hD\rangle^\ell u \in L^r \iff u\in W^{-k,r} W^{\ell,r }.\]
\end{subsection}

\begin{subsection}{Tangential Calculus}
We denote the H\"ormander symbols by $S^{\ell}_1 (\R^n)$. We also consider symbols in the class $S^k_0(\R^n)$. In this article we will work with product symbols of the form $ba(x',\xi) \in S^k_1(\R^{n-1})S^\ell_j(\R^n) := S^k_1 S^\ell_j$ where $b(x',\xi')\in S^k_{1}(\R^{n-1})$ and $a(x',\xi)\in S^\ell_j(\R^n)$ for $j =0,1$. Observe that if $a(x',\xi) \in S^{k}_1 S^\ell_j$, then derivatives with respect to either $x'$ or $\xi$ are a finite sum of symbols in $S^{k}_1 S^\ell_j$:
\begin{eqnarray}
\label{differentiating product symbols}
\partial_{x}^\alpha : S^{k}_1 S^\ell_j \to {\mathrm{span}}(S^{k}_1 S^\ell_j)\ \ \, \partial_{\xi}^\alpha : S^{k}_1 S^\ell_j \to {\mathrm{span}}(S^{k}_1 S^\ell_j).
\end{eqnarray}

We begin with the following Calder\'on-Vaillancourt type estimate for (classical) $\Psi$DO with symbols in $S_1^0(\R^n)$ which can be obtained by following the argument of Theorem 9.7 in \cite{wong}.
\begin{prop}
\label{calderonvaillancourt}
Let $a(x,\xi)$ be a symbol in $S_1^0(\R^n)$. Then for all $1<r<\infty$ 
\begin{eqnarray}
\label{classical estimate}
\|a(x,D) u\|_{L^r} \leq C_{r,n} \sum\limits_{|\alpha|\leq k(n), |\beta| \leq k(n)} p_{\alpha,\beta}(a) \|u\|_{L^r}
\end{eqnarray}
where $p_{\alpha,\beta}$ is the semi-norm defined by $p_{\alpha,\beta}(a) := \sup\limits_{x,\xi} |\partial_x^\alpha \partial_\xi^\beta a(x,\xi)| \langle \xi\rangle^{|\beta|}$ and $k(n)\in \N$ depends on the dimension only.
\end{prop}
We shall henceforth denote by $k(n)$ to be the smallest integer for which Proposition \ref{calderonvaillancourt} holds.
Note that in $\R^n$ there is a relation between classical and semiclassical quantization of a symbol $a\in S^\infty$ given by
\[ Op_h(a)u (x) = \sqrt{h}^{-n} A_h u_h (x/\sqrt{h})\]
where $u_h$ is defined by $({\mathcal F} u_h)(\xi) = ({\mathcal F}u) (\xi/\sqrt{h})$ and $A_h = a_h(x,D)$ for $a_h(x,\xi) := a(\sqrt{h} x,\sqrt{h}\xi)$ ($\mathcal F$ denotes the {\bf classical} Fourier transform). This identity combined with estimate \eqref{classical estimate} gives us a semiclassical version of Calder\'{o}n-Vaillancourt: for all $1<r<\infty$ and $h>0$ sufficiently small,
\begin{eqnarray}\label{semiclassical calderonvaillancourt}\|Op_h (a) u\|_{L^r} \leq \sum\limits_{|\alpha|,|\beta| \leq k(n)} p_{\alpha,\beta}(a) \|u\|_{L^r} + C\sqrt{h}\|u\|_{L^r}.\end{eqnarray}

  For symbols in $S^{k}_1 S^{-\ell}_1 \cup S^{k}_1 S^{-k(n) -\ell}_0$, we have the following mapping properties.
\begin{prop}
\label{sobolev mapping}
If $b(x',\xi')\in S^{k}_1$ and $a(x',\xi) \in S^{\ell}_1 \cup S^{-k(n)+\ell}_0$ then \[ba(x', hD) : W^{m,r} W^{l,r} \to W^{m-k,r} W^{l - \ell,r}\] with norm
{\tiny \[\|ab(x', hD)\| \leq C \sup\limits_{\substack{z', \xi,\\ |\alpha| \leq k(n)}} |(1+ \Delta_{z'})^N \partial_\xi^\alpha a(z', \xi)| \langle \xi\rangle ^{|\alpha|-\ell}\sup\limits_{\substack{z', \xi,\\ |\alpha| \leq k(n)}} |(1+ \Delta_{z'})^N \partial_{\xi'}^\alpha b(z', \xi')| \langle \xi'\rangle ^{|\alpha|-k}.\]}

\end{prop}
%\end{subsection}

In addition, we have the following compositional calculus result.
\begin{prop}
\label{sobolev composition}
If $a\in S^{k_1}_1S_1^{\ell_1} \cup S^{k_1}_1 S^{-k(n) +\ell_1}_0$ and $b\in S^{k_2}_1S_1^{\ell_2} \cup S^{k_2}_1 S^{-k(n) +\ell_2}_0$ then
\[  b(x'hD)a(x',hD) = ab(x',hD) + h\sum\limits_{|\alpha|=1} (\partial_{\xi}^\alpha b\partial_{x'}^\alpha a)(x',hD) + h^2 m(x',hD)\]
where $m(x',hD): W^{k,r}W^{\ell,r} \to W^{k - k_1-k_2,r}W^{\ell-\ell_1-\ell_2,r}$. 
\end{prop}

For proofs of Proposition \ref{sobolev mapping} and Proposition \ref{sobolev composition}, see the Appendix.

%\lt{put in a remark to cover ourselves for weighted $L^2$ spaces}
\begin{remark} We have omitted stating the mapping properties on $H^k_\delta$ spaces since $S^k_0 S^\ell_0 \subset S_0^{k+\ell}(\R^n)$ and the calculus for these symbols on weighted $L^2$-Sobolev spaces are well documented. See for example \cite{salothesis}.
\end{remark}
\end{subsection}
\end{section}

\begin{section}{Heat Flow}

Define coordinates on $\R^n$ and let $\R^n_+$ denote the upper half space $\{x_n > 0\}$. Let $F(x',\xi') \in S^1_1(\R^{n-1})$, and define the semiclassical pseudodifferential operator 
\begin{equation}\label{J}
j(x',hD) = h\partial_{x_n} + F(x',hD')
\end{equation}
on $\R^n$. It follows by considering the $\xi'$ and $\xi_n$ direction separately and applying the semiclassical Calder\'{o}n-Vaillancourt theorem that $j(x',hD)$ is a bounded operator $j(x',hD):W^{1,r}(\R^n) \rightarrow L^r(\R^n)$ for $1 < r < \infty$. As we will see in the following section, one of the factors of the conjugated Laplacian has this form. In this section we will prove some basic facts about the existence and $L^p$ mapping properties of the inverse of such an operator. This extends the $L^2$ theory explained in \cite{ChuND}. 

To obtain an inverse, we will assume that $F$ obeys the ellipticity condition
\begin{equation}\label{EllipticF}
c\langle \xi' \rangle  \leq \mathrm{Re}F(x',\xi') \leq C\langle \xi' \rangle
\end{equation}
uniformly in $x'$, for some constants $c, C > 0$.  This ensures that the principal symbol
\[
j(x,\xi) := i\xi_n + F(x',\xi') 
\]
is uniformly elliptic. We will also assume a finiteness condition on $F$: that there exists $X' > 0$ such that for $|x'| > X',$
\begin{equation}\label{FiniteF}
\grad_{x'} F(x',\xi') = 0.
\end{equation}
We need an extra condition to ensure that the symbol $j^{-1}$ is in the suitable calculus. We assume that there exists a first order symbol $i\xi_n + F_-(x',\xi')$ with compact characteristic set, such that $D_{x'} F_-(x',\xi')$ is supported in $|x'| <X'$, and 
\[
(i\xi_n + F) (i\xi_n + F_-) = p(x',\xi) + a_0
\]
where $p(x',\xi)$ is a second order polynomial in $\xi$ with compact characteristic set and $a_0 \in S^{-\infty}(\R^{n-1})$. 

The reason why we need this extra assumption is that $(i\xi_n + F)^{-1}$ is not in general in the class $S^{-1}_1(\R^n)$. However, if $\chi\in C^\infty_0(\R^{n})$ is identically $1$ on a neighbourhood containing the characteristic sets of $i\xi_n + F_-$ and $p$, then we can derive the following expansion:
\begin{eqnarray*}
(1- \chi(\xi))(i\xi_n + F)^{-1} = (1-\chi(\xi))(i\xi_n + F_-)\left(\frac{1}{p(x',\xi)} - \frac{a_0}{p(p+a_0)}\right).
\end{eqnarray*}
Since $\chi$ is identically one on the characteristic set of $p$, it follows that $(1 - \chi(\xi))/p(x',\xi)$ is a symbol in $S^{-2}_1(\R^n)$, and so
\[
(1- \chi(\xi))(i\xi_n + F)^{-1} = (i\xi_n + F_-)\left(S^{-2}_1 - \frac{a_0(1 - \chi)}{p(p+a_0)}\right).
\] 
Now observe that $\frac{a_0}{p(p+a_0)} = \frac{a_0}{p^2} - \frac{a_0^2}{p^2(p+a_0)}$, and we can repeat this procedure indefinitely to obtain
\begin{equation*}
(1- \chi(\xi))(i\xi_n + F)^{-1} = (i\xi_n + F_-)\Big(S^{-2}_1 + a_0 S_1^{-4} +\dots+ a_0^m S^{-k(n) -1}_1 + a_0^{m+1} S_0^{-k(n) -2}\Big)
\end{equation*}
where we are using $S^k_j$ to represent a symbol from the class $S^k_j(\R^n)$. 
Now $(i\xi_n + F_-)S^{-2}_1 \in S^{-1}_1 + S^1_1S^{-2}_1$, and the same holds for $(i\xi_n + F_-)(a_0 S_1^{-4} +\dots+ a_0^m S^{-k(n) -1}_1)$. Finally, $(i\xi_n + F_-)a_0^{m+1} S_0^{-k(n) -2} \in {\mathrm{span}}(S^{-\infty} S_{0}^{-1-k(n)} )$, so     
\begin{equation}\label{expansion of inverse}
(1- \chi(\xi))(i\xi_n + F)^{-1} \in {\mathrm{span}}( S^0_1 S^{-1}_1 + S^{-\infty} S_{0}^{-1-k(n)} + S_1^1 S^{-2}_1).
\end{equation}

Meanwhile $\chi(\xi) j^{-1}(x',\xi)\in S^{-\infty}(\R^n)$, so we can use \eqref{expansion of inverse} in conjunction with Proposition \ref{sobolev mapping} to get that 
%\lt{added in the statement about weighted spaces}
\begin{eqnarray}
\label{j parametrix}
j^{-1}(x',hD) : L^2_\delta \to H^1_\delta,\ \ \delta\in \R,\ \ j^{-1}(x',hD) : L^r \to W^{1,r},\ \ 1<r<\infty.
\end{eqnarray}

The operator $j^{-1}(x',hD)$ also turns out to have desirable support properties.  

\begin{lemma}
\label{support of j^-1}
If $u\in L^r(\R^n)$ is supported only in $\{x_n \geq 0\}$ then $j^{-1}(x',hD)u \in W^{1,r}(\R^n)$ has trace zero along $\{x_n = 0\}$ and vanishes identically on the set $\{x_n <0\}$.
\end{lemma}

\begin{proof}
For $u \in C^{\infty}_c(\R^n)$, we can write %\lt{I have changed the notation here to ${\mathcal F}^h$ to denote the semiclassical transform}
\[
j^{-1}(x',hD)u(x) = h^{-n}\int_{\R^{n-1}} \int_{\R}\frac{{\mathcal F}^h{u}(\xi)}{i\xi_n + F(x',\xi')} e^{\frac{i}{h}x\cdot \xi} d\xi_n \, d\xi'
\]  
where ${\mathcal F}^h$ is the semiclassical Fourier transform. We can write out the Fourier transform in the $x_n$ variable to get
\[
j^{-1}(x',hD)u(x) = h^{-n}\int_{\R^{n-1}} \int_{\R} \int_{\R}\frac{\mathcal{F}^h_{x'}u(\xi',t)}{i\xi_n + F(x',\xi')} e^{\frac{i}{h}(x_n-t)\xi_n} d\xi_n \, dt \, e^{\frac{i}{h}x'\cdot \xi'} d\xi'.
\]
Now we can use the residue theorem to evaluate the $d\xi_n$ integral explicitly, and we get
\[
j^{-1}(x',hD)u(x) = h^{-n}\int_{\R^{n-1}} \int_{-\infty}^{x_n}\mathcal{F}^h_{x'}u(\xi',t)e^{\frac{t-x_n}{h} F(x',\xi')} dt \, e^{\frac{i}{h}x'\cdot \xi'} d\xi'.
\]
For $u \in C^{\infty}_c(\R^n)$, the lemma follows immediately from this representation.  Then the lemma holds for general $u \in L^r(\R^n)$ by using a density argument involving the bounds in \eqref{j parametrix}.  

\end{proof}
Henceforth we will refer to the support property given in Lemma \ref{support of j^-1} as ``preserving support in the $x_n$ direction". 

We can turn $j^{-1}(x',hD)$ into a proper inverse. We first prove a composition type lemma for the operator $j^{-1}(x',hD)$.
\begin{lemma}
\label{composing with j^-1}
Let $a(x',\xi') \in S^{1}_1(\R^{n-1})$. Then 
\[a(x',hD') j^{-1}(x',hD) = (aj^{-1})(x',hD) + h\sum\limits_{|\alpha|=1} (j^{-2} \partial_{\xi'}^\alpha a \partial_{x'}^\alpha F )(x',hD)+ h^2m(x',hD)\] where $m(x',hD)$ and $\sum\limits_{|\alpha|=1} (j^{-2} \partial_{\xi'}^\alpha a \partial_{x'}^\alpha F)(x',hD)$ map $L^{r} \to L^r$ with norm bounded by a constant independent of $h$. Furthermore, the commutator $[a(x',hD'), j^{-1}(x',hD)]  = hm(x,hD)$ with \[m(x,hD): L^r\to L^r, \ \ m(x,hD): L^2_\delta\to L^2_\delta.\]
\end{lemma}

\begin{proof}
The expansion \eqref{expansion of inverse} allows us to write $j^{-1}(x',\xi)$ as span of elements in 
\[  S^0_1S^{-1}_1 + S^{-\infty} S_{0}^{-1-k(n)} + S_1^1 S^{-2}_1 .\]
We can therefore apply Proposition \ref{sobolev composition} to each term to obtain 
\begin{eqnarray*}
%\label{approximate inverse for J}
a(x',hD') j^{-1}(x',hD) = aj^{-1}(x',hD) + h m_1(x',hD) + h^2 m_2(x',hD)
\end{eqnarray*}
where \[m_1(x',\xi) = \sum\limits_{|\alpha|=1} (j^{-2}\partial_{\xi'}^\alpha a \partial_{x'}^\alpha F)(x',\xi), \, \mbox{ and }  \ \  m_2(x',hD) : L^r \to L^r.\] 
Using expansion \eqref{expansion of inverse} again we see that $m_1(x',\xi)$ is a symbol in the span of \[S^1_1S^{-1}_1+S^{-\infty}S_0^{-k(n)-1}+S_1^2 S_1^{-2}.\] Therefore, it maps $L^r\to L^r$ by Proposition \ref{sobolev mapping} and the fact that $W^{-k,r} W^{\ell,r } \subset W^{l-k, r}(\R^n)$.
To obtain the commutator statement, repeat the argument for the composition $ j^{-1}(x',hD)a(x',hD')$.
\end{proof}

Now we can use $j^{-1}$ to build a proper inverse which preserves support in the $x_n$ direction.  Moreover, the inversion can still be carried out even if $j$ is perturbed by a small tangential operator $hF_0$.  

\begin{prop} %\lt{ I see that $hF_0$ is introduced here. I am OK with this way of presenting the operator}
\label{J inverse}
Suppose $F_0(x',\xi') \in S^{0}_1(\R^{n-1})$ obeys the same finiteness condition \eqref{FiniteF} as $F$, and consider the operator 
\[
J := j(x',hD) + hF_0(x',hD). 
\]
For $h>0$ sufficiently small there exists an inverse  $J^{-1}: L^r \to W^{1,r}$ of the form
\[J^{-1}=j^{-1}(x', hD)(1 + h m_1(x',hD) + h^2 m_2(x',hD))^{-1}\] 
where $m_1(x',hD),m_2(x'hD): L^r \rightarrow L^r$. %$m_1(x',hD) := \sum\limits_{|\alpha|=1} (j^{-2}\partial_{\xi'}^\alpha F \partial_{x'}^\alpha F)(x',hD) + (j^{-1}F_0)(x',hD): L^r\to L^r$ and $m_2(x',hD) : L^r \to L^r$. 

Furthermore, $J^{-1}$ preserves support in the $x_n$ direction: if the support of $u\in L^r$ is contained in $x_n \geq 0$ then $J^{-1} u$ has vanishing trace on $\{x_n = 0\}$ and vanishes identically when $\{x_n <0\}$. The same holds for mapping properties on $H^k_\delta$ spaces.
\end{prop}
\begin{proof} We write 
\[
Jj^{-1}(x',hD) = h\partial_n j^{-1}(x',hD) + F(x',hD') j^{-1}(x',hD) + hF_0(x',hD') j^{-1}(x',hD).
\] 
We can apply Proposition \ref{sobolev composition} to the first term, using the expansion \eqref{expansion of inverse} for $j^{-1}$, and Lemma \ref{composing with j^-1} to the second and third terms to obtain  
\begin{eqnarray}
\label{approximate inverse for J}
J j^{-1}(x',hD) = 1 + h m_1(x',hD) + h^2 m_2(x',hD)
\end{eqnarray}
where 
\begin{equation}\label{JInverseRemainder}
m_1(x',\xi) = \sum\limits_{|\alpha|=1} (j^{-2}\partial_{\xi'}^\alpha F \partial_{x'}^\alpha F)(x',\xi) + (j^{-1}F_0)(x',\xi'),
\end{equation}
and $m_2(x',hD) : L^r \to L^r.$ Using expansion \eqref{expansion of inverse} again we see that $m_1(x',\xi)$ is a symbol in the span of \[S^1_1S^{-1}_1+S^{-\infty}S_0^{-k(n)-1}+S_1^2 S_1^{-2}.\] Therefore, it maps $L^r\to L^r$ by Proposition \ref{sobolev mapping} and the fact that $W^{-k,r} W^{\ell,r } \subset W^{l-k, r}(\R^n)$.

Observe that in equation \eqref{approximate inverse for J} since $J$ is a differential operator in the $x_n$ direction, it preserves support in the $x_n$ direction when acting on $W^{1,r}$. The operator $j^{-1}(x',hD) :L^r\to W^{1,r}$ preserves support in $x_n$ by Lemma \ref{support of j^-1} and thus the left side preserves support in the $x_n$ direction. We may conclude from this that the right side preserves $x_n$ support as well and in particular $h m_1(x',hD) + h^2 m_2(x',hD)$ preserves $x_n$ support. This means that inverting the right-side by Neumann series preserves support in the $x_n$ direction.
\end{proof}

One final consequence of the structure of $J^{-1}$ we obtained in Proposition \ref{J inverse} is the following disjoint support property:
\begin{lemma}
\label{disjoint support for J inverse}
Let ${\bf 1}_{\R^n_-}$ be the indicator function for $x_n \leq 0$ and $\epsilon >0$. Then for all $f\in L^r(\R^n)$,
\[
\|J^{-1}{\bf 1}_{\R^n_-}f \|_{W^{1,r}(\{x_n \geq \epsilon\})} \leq C_\epsilon h^2 \|f\|_{L^r}.
\]
\end{lemma}

\begin{proof}
Let $\zeta_{\epsilon}(x_n)$ be a smooth cutoff function which is identically one on $\{x_n \geq \epsilon \}$ and identically zero on an open set containing $\{x_n \leq 0\}$. Then 
\[
\|J^{-1}{\bf 1}_{\R^n_-}f \|_{W^{1,r}(\{x_n \geq \epsilon\})} \leq \|\zeta_{\epsilon}J^{-1}{\bf 1}_{\R^n_-}f \|_{W^{1,r}(\R^n)}.
\]
Therefore it suffices to show that 
\[
\|\zeta_{\epsilon}J^{-1}{\bf 1}_{\R^n_-}f \|_{W^{1,r}(\R^n)} \leq C_\epsilon h^2 \|f\|_{L^r}.
\]
From Proposition \ref{J inverse}, we have that 
\[
J^{-1} = j^{-1}(x' hD)(1 + hm_1(x',hD) + h^2m_2(x',hD))^{-1}
\]
where $(1 + hm_1(x',hD) + h^2m_2(x',hD))^{-1}$ is given by the Neumann series
\[
(1 + hm_1(x',hD) + h^2m_2(x',hD))^{-1} = 1 + \sum_{k=1}^\infty (hm_1(x',hD) + h^2m_2(x',hD))^k.
\]
Therefore, by \eqref{j parametrix} we can write 
\[
J^{-1} = j^{-1}(x', hD)(1 + hm_1(x',hD)) + h^2M(x',hD)
\]
where $M: L^r \rightarrow W^{1,r}$ is bounded uniformly in $h$.  Using this expression for $J^{-1}$ it suffices to show that
\[\zeta_\epsilon j^{-1}(x', hD)(1 + hm_1(x',hD)) {\bf 1}_{\R^n_-}: L^r(\R^n)\to W^{1,r}(\R^n)\]
with norm bounded by $O(h^2)$.
%\[
%\zeta_\epsilon J^{-1}{\bf 1}_{\R^n_-}f = \zeta_\epsilon  j^{-1}(x', hD)(1 + hm_1(x',hD)){\bf 1}_{\R^n_-}f.
%\]
We will demonstrate this for the principal part $\zeta_\epsilon j^{-1}(x',hD) {\bf 1}_{\R^n_-}$ and leave the lower order term, which can be written out explicitly using \eqref{JInverseRemainder}, to the reader. By using \eqref{expansion of inverse} we see that the symbol 
\[j^{-1} \in {\mathrm{span}} (S^0_1S^{-1}_1 + S^{-\infty}S^{-1-k(n)}_0 + S^1_1S^{-2}_1).\]
We will only show the estimate for $ \zeta_\epsilon Op_h(S^1_1S^{-2}_1)  {\bf 1}_{\R^n_-}$ and the others are treated in the same way.
Suppose $b \in S^1_1(\R^{n-1})$ and $a \in S^{-2}_1(\R^n)$, by Proposition \ref{sobolev composition} we see that
{\Small\[ \zeta_\epsilon ba(x',hD) {\bf 1}_{\R^n_-} = \zeta_\epsilon b(x',hD')a(x',hD) {\bf 1}_{\R^n_-} + h\zeta_\epsilon\sum_{|\alpha| = 1} (\partial_{\xi'} b)(x',hD') (\partial_{x'} a)(x', hD) {\bf 1}_{\R^n_-} + h^2 m(x, hD)\]}
where $m(x',hD) : L^r \to W^{-1,r}W^{2,r} \subset W^{1,r}(\R^n)$ by \eqref{mixed space embedding}.

Since $\zeta_{\epsilon}$ is a function of $x_n$ only, it commutes with operators from $S^k_1(\R^{n-1})$, and thus estimating $\zeta_\epsilon ba(x',hD) {\bf 1}_{\R^n_-} $ with $b \in S^1_1(\R^{n-1})$ and $a \in S^{-2}_1(\R^n)$ amounts to estimating terms of the form
$ \zeta_\epsilon Op_h(S^{-2}_1(\R^n)) {\bf 1}_{\R^n_-}$. Standard disjoint support properties of $\Psi$DO then give the desired estimates.
\end{proof}

\end{section}

\begin{section}{Green's Functions on $\R^n$}

The purpose of this discussion is to find a way to invert 
\[h^2\Delta_\phi  := h^2 e^{-\phi/h}\Delta e^{\phi/h},\ \ \phi(x) := x_n\]
with a suitable boundary condition and good $L^{p'}\to L^p$ estimates. We begin with the operator on $\R^n$ given by the Fourier multiplier $\frac{1}{|\xi|^2 + 2i\xi_n -1}$. We give a semiclassical formulation of an estimate established in Sylvester-Uhlmann \cite{su}.
\begin{lemma}%\lt{realized that the map is from $L^2_\delta\to H^2_{\delta-1}$. Made changes in this section to reflect this fact}
\label{sylvester uhlmann}
The Fourier multiplier $\frac{1}{|\xi|^2 + 2i\xi_n -1}$ maps $L^2_\delta \to H^2_{\delta-1}$ for $\delta >0$ with norm bounded by $h^{-1}$.
\end{lemma}
\begin{proof}
Consider the multiplier given by $\frac{1}{|\xi|^2 + 2i\xi_n + 2\xi_1}$. By the result of \cite{su}, 
\[Op_h\left(\frac{1}{|\xi|^2 + 2i\xi_n + 2\xi_1}\right) : L^2_\delta \to_{h^{-1}} H^2_{\delta-1}.\]
Observe that $|\xi|^2 + 2i\xi_n + 2\xi_1= |\xi_1 + 1|^2 + \sum\limits_{j = 2}^n |\xi_j|^2 + 2i\xi_n -1$. Since shifting in the Fourier coordinate is equivalent to multiplying by a complex linear phase,
\[Op_h\left(\frac{1}{|\xi|^2 + 2i\xi_n -1}\right) = e^{-ix_1/h} Op_h\left(\frac{1}{|\xi|^2 + 2i\xi_n + 2\xi_1}\right) e^{ix_1/h}\]
and the proof is complete.
\end{proof}

It turns out that the Fourier multiplier $\frac{1}{|\xi|^2 + 2i \xi_n -1}$ also satisfies $L^{p'}\to L^p$ estimates for $p = \frac{2n}{n-2}$ and $p' = \frac{2n}{n+2}$. We describe below the semiclassical formulation of a result by Kenig-Ruiz-Sogge \cite{krs} and Chanillo \cite{ch} -- see also Haberman \cite{Hab}.
\begin{lemma}
\label{krs estimate}
The Fourier multiplier satisfies the estimate
\[\left\|Op_h\left(\frac{1}{|\xi|^2 + 2i \xi_n -1}\right)u\right\|_{L^{p}} \leq \frac{C}{h^2} \|u\|_{L^{p'}}\]
for all $u\in L^{p'}(\R^n)$.
\end{lemma}
\begin{proof}
We begin with a classical estimate for $Op_1\left(\frac{1}{|\xi|^2 + 2i \xi_n -1}\right)$ due to \cite{krs}. Let $u\in\S$ be a Schwartz function satisfying $\hat u(\xi',\xi_n) = 0$ for whenever $\xi_n$ is close to zero. For these $u$, we have $Op_1\left(\frac{1}{|\xi|^2 + 2i \xi_n -1}\right)u \in W^{2,p'}(\R^n)$ and we can therefore apply Theorem 2.4 of \cite{krs} to obtain
\begin{eqnarray}
\label{krs estimate on test func}
\left\|Op_1\left(\frac{1}{|\xi|^2 + 2i \xi_n -1}\right)u \right\|_{L^{p}} \leq C\|u\|_{L^{p'}}.\end{eqnarray}
We would like to use a density argument to show that the above holds for all $u\in L^{p'}$. Indeed, let $u\in\S$ be any Schwartz function and define for all $\delta>0$ the Schwartz function 
\[
\hat u_\delta (\xi) := \theta_\delta(\xi_n)\hat u := \theta(\xi_n/\delta)\hat{u}
\]
where $\theta : \R \to [0,1]$ is a smooth bump function which is identically $1$ near the origin. By the dominated convergence theorem and Plancherel one sees that $\lim\limits_{\delta\to 0} \|u_\delta\|_{L^2} = 0$. For the $L^1$ norm, observe that 
\[u_\delta (x) = \int_{\R} e^{i\xi_n x_n} \theta_\delta(\xi_n) \int_{\R} e^{-i\xi_n t} u(x',t)dt d\xi_n = \int_{\R} \check\theta_\delta (s) u(x',x_n-s) ds\]
so one has $\|u_\delta\|_{L^1} \leq \delta \|u\|_{L^1} \int_{\R}|\check\theta (\delta s)|ds \leq C\|u\|_{L^1}$. Riesz-Thorin interpolation then yields that $\|u_\delta\|_{L^{p'}} \to 0$ for $p' = 2n/(n+2)$. The function $u - u_\delta$ is then an element of $\S$ whose Fourier transform vanishes in a neighbourhood of $\xi_n = 0$ which converges to $u$ in $L^{p'}$ and thus \eqref{krs estimate on test func} is valid for all $u \in L^{p'}$ by density.

Denote by $u_h(x) := u(hx)$ and insert $u_h$ into the estimate \eqref{krs estimate on test func} in place of $u$. We get
\[
\left\|\Big(Op_h(\frac{1}{|\xi|^2 + 2i \xi_n -1})u\Big)_h \right\|_{L^{p}} \leq C\|u_h\|_{L^{p'}}.
\]
Making the change of variable $y = hx$ and computing the norms on both sides we get the desired semiclassical estimate stated in the Lemma. The $h^{-2}$ factor arises from the fact that $1/p' - 1/p = 2/n$. 
\end{proof}

%\begin{remark} The authors would like to thank Professor Chanillo for showing them this argument.\end{remark}

In order to deal with domains with non-flat boundaries, we will actually need to deal with domains ``flattened'' by a coordinate change of the type 
\begin{eqnarray}
\label{change of var}
\gamma : (y_n,y') \mapsto (x_n , x') = (y_n - f(y'), y').\end{eqnarray}
Under this change of variables, the conjugated Laplacian $-e^{-y_n/h} (\sum_{j=1}^n h^2\partial_{y_j}^2) e^{y_n/h}$ becomes the operator
\[
h^2 \tilde\Delta_\phi  =  Op_h((1 + |K|^2)\xi_n^2 -2\xi_n(i-\xi'\cdot K) -(1- |\xi'|^2) )
\] 
where $K(x') := \nabla f (x')$. The next proposition concerns the Green's function for $h^2 \tilde\Delta_\phi$.
 \begin{prop}
 \label{green with no boundary}
 The Green's function defined by $\tilde G_\phi := \gamma^* G_\phi$ satisfies $h^2\tilde\Delta_\phi \tilde G_\phi = Id$ and has the bounds \[\|\tilde G_\phi\|_{L_{\delta}^2 \to H_{\delta-1}^2} \leq Ch^{-1},\ \ \|\tilde G_\phi\|_{L^{p'} \to L^{p}} \leq Ch^{-2}.\]
 Furthermore, we can split $\tilde G_\phi = \tilde G_\phi^c + (\tilde G_\phi - \tilde G_\phi^c)$ such that $(\tilde G_\phi - \tilde G_\phi^c)$ is a $\Psi$DO with symbol in $S^{-2}_1(\R^n)$ and 
 \[\|\tilde G_\phi^c\|_{L^{p'} \to W^{k,p}} \leq Ch^{-2},\ \ \|\tilde G^c_\phi\|_{L_{\delta}^2 \to H_{\delta-1}^k} \leq Ch^{-1} \ \forall k \in \N.\]
  \end{prop}
\begin{proof}
By Lemmas \ref{krs estimate} and \ref{sylvester uhlmann} the operator $G_\phi$ defined by the semiclassical multiplier $\frac{1}{|\xi|^2 - 2i\xi_n -1}$ satisfies $h^2 \Delta_\phi G_\phi = I$, $G_\phi : L^2_\delta(\R^n) \to H^2_{\delta-1}(\R^n)$ with $O(h^{-1})$ norm and $L^{p'}(\R^n) \to L^p(\R^n)$ with $O(h^{-2})$ norm. 

The multiplier of $G_\phi$ is constant coefficient so one can write $G_\phi = G^c_\phi + (G_\phi - G^c_\phi)$ where \[G^c_\phi  = G_\phi \chi_0(hD) = \chi_1(hD)G_\phi \chi_0(hD)\] where $\chi_0, \chi_1 \in C^{\infty}_c(\R^n)$, with $\chi_0$ identically $1$ in the ball of radius $2$ and $\chi_1$ identically $1$ on the support of $\chi_0$.

Since the characteristic set of $G_\phi$ is disjoint from the support of $1-\chi_0$, the operator $ (G_\phi - G^c_\phi): L^r(\R^n) \to W^{2,r}(\R^n)$ is a $\Psi$DO with symbol in $S^{-2}_1(\R^n)$.

The mapping properties of $G^c_\phi$ come from the mapping properties of $G_\phi$ and the fact that $\chi_1(hD)$ has compactly supported symbol. 

The estimates for the pull-back operator $\tilde G_\phi$ follows naturally from the estimates for $G_\phi$ since the Jacobian of $\gamma$ is identity outside of a compact set.
\end{proof}

The characteristic set of $G_{\phi}$ lies in the sphere $|\xi'| = 1$, and so in particular if $G^c_{\phi}$ is multiplied by a Fourier side cutoff function supported away from that sphere, the resulting operator is well behaved.  The following lemma makes this somewhat more precise.

\begin{lemma}
\label{disjoint from char}
Let $\t \rho(\xi')$ be a smooth function with support compactly contained in $|\xi'| <1$. Then $\t \rho \t G_\phi^c = Op_h(S^{-\infty}(\R^n)) + h m(x',hD)\t G_\phi$ for some $m(x',\xi)\in S^{-\infty}(\R^n)$.
\end{lemma}
\begin{proof}
By Proposition \ref{green with no boundary} $\t G_\phi^c = \gamma^*( \chi G_\phi)$ where $\gamma^*$ is the pull-back by the diffeomorphism given by $(x',x_n) \mapsto ( x', x_n - f(x'))$. We compute 
{\small\[\t \rho(hD') \gamma^*(\chi(hD) G_\phi) = \t \rho(hD') \gamma^*(\chi(hD))\t G_\phi=( \t \rho(hD')\t \chi(x',hD) + hOp_h(S^{-\infty}(\R^n)))\t G_\phi\]}
where $\t \chi (x',\xi) = \chi (D\gamma(x')^T \xi)$ is the pull-back symbol. By the composition formula in Proposition \ref{sobolev composition},
{\Small \[\t \rho(hD') \gamma^*(\chi(hD) G_\phi) = (Op_h(\t \rho \t \chi) + hOp_h(S^{-\infty}(\R^n)))\t G_\phi = \gamma^*(Op_h(\chi\rho)G_\phi) + hOp_h(S^{-\infty}(\R^n))\t G_\phi\]}
where $\rho(x',\xi) = \t \rho (\xi' + \xi_n K(x'))$ is the push-forward symbol. Observe that since $G_\phi$ is constant coefficient, $Op_h(\chi\rho)G_\phi = Op_h\left(\frac{\chi(\xi)\rho(x',\xi)}{|\xi'|^2 + \xi_n^2 -1 + 2i\xi_n}\right)$. Since $\t\rho(\xi')$ vanishes in an open neighbourhood of $|\xi'|=1$, the symbol \[\frac{\chi(\xi)\rho(x',\xi)}{|\xi'|^2 + \xi_n^2 -1 + 2i\xi_n}\]belongs to $S^{-\infty}(\R^n)$.
\end{proof}

\begin{subsection}{Modified Factorization}
To add boundary determination to the Green's function, we want to take advantage of the fact that $h^2\tilde{\Delta}_{\phi}$ factors into two parts, one of which is elliptic and resembles the operator described in Section 3.  

Indeed, the symbol of $\frac{1}{1+K^2}h^2\tilde \Delta_\phi$ factors formally as
{\tiny\begin{eqnarray*}\xi_n^2 -2\xi_n\frac{(i-\xi'\cdot K)}{ 1+ |K|^2} -\frac{(1- |\xi'|^2)}{1+|K|^2}&=& \left(\xi_n - i\left( \frac{(1 + iK\cdot \xi')-\sqrt{(1 +i K\cdot \xi')^2 - (1-|\xi'|^2)(1+ |K|^2)}}{1+|K|^2}\right)\right)\times\\&&\left(\xi_n - i\left( \frac{(1 + iK\cdot \xi')+\sqrt{(1 +i K\cdot \xi')^2 - (1-|\xi'|^2)(1+ |K|^2)}}{1+|K|^2}\right)\right) \end{eqnarray*}} 
and the second factor here is elliptic. The problem is that the square root is not smooth at its branch cut, so this does not give a proper factorization at the operator level. The obvious thing to do is to take a smooth approximation to the square root, but for our purposes we will require something more subtle.

We take the branch of the square root that has non-negative real part, and seek to avoid the branch cut, which happens when the argument of the square root lies on the negative real axis. From examination of the square root, we see that this occurs when $K\cdot \xi '= 0$ and $|\xi'|^2 \leq |K|^2(1 + |K|^2)^{-1}$. By ensuring that $\xi'$ avoids this set, we can guarantee that the argument of the square root stays away from the branch cut. % Calculating more carefully gives us the following quantitative version of that statement.

Thus let $0<c< c'<1$ be a constant such that $\frac{|K|^2}{1+|K|^2} <c$ for all $x'$ and let $\tilde\rho_0(\xi')$ be a smooth function in $\xi'$ such that $\tilde \rho_0 = 1$ for $|\xi'|^2 \leq c$ and ${\supp} (\tilde\rho_0)  \subset\subset B_{\sqrt{c'}}$. Introduce a second cutoff $\tilde \rho$ such that it is identically $1$ on $|\xi'|^2 \leq c'$ but ${\supp} (\tilde\rho)  \subset\subset B_{1}$. Observe that
\begin{eqnarray}
\label{supp of rho tilde}
\inf_{\xi\in {\supp}{\t\rho},\ x'\in \R^{n-1}} \left|\xi_n^2 -2\xi_n\frac{(i-\xi'\cdot K)}{ 1+ |K|^2} -\frac{(1- |\xi'|^2)}{1+|K|^2}\right| >0.
\end{eqnarray}

%Lemma \ref{1-rho avoids branch cut} shows
Since the branch cut of the square root occurs when $|\xi'|^2 \leq |K|^2(1 + |K|^2)^{-1}$, it follows that for $\xi'$ in the support of $1-\tilde\rho_0$, the function $(1 + iK\cdot \xi')^2 - (1-|\xi'|^2)(1+ K^2)$ stays uniformly away from the branch cut of the square root. As such we may define
\begin{eqnarray}
\label{square root}
r := (1- \t\rho_0) \sqrt{(1+iK\cdot \xi' )^2 -(1-|\xi'|^2)(1+ |K|^2)}
\end{eqnarray}
and factor
{\Tiny \begin{eqnarray}
\label{decomposition}\nonumber
\xi_n^2 -2\xi_n\frac{(i-\xi'\cdot K)}{ 1+ |K|^2} -\frac{(1- |\xi'|^2)}{1+|K|^2}= {(\xi_n - \t a_- +h m_0) }{(\xi_n - \t a_+ -h m_0)} + \t a_0 + h\sum\limits_{|\alpha| = 1} \partial_{\xi'}^\alpha \t a_-\partial_{x'}^\alpha \t a_+ - h{m_0 \t a_-}+ h^2 m_0^2\\
\end{eqnarray}}
with $m_0(x',\xi') := -\t a_+^{-1}\sum\limits_{|\alpha| = 1} \partial_{\xi'}^\alpha \t a_-\partial_{x'}^\alpha \t a_+$.
Here the $\t a_\pm$ and $\t a_0$ are defined by 
{\Small\begin{eqnarray}
\label{the a's}
\t a_\pm =  i\left( \frac{(1 + iK\cdot \xi')\pm r}{1+|K|^2}\right)\ \ {\rm and}\ \ \t a_0 =  \frac{(1 +i K\cdot \xi')^2 - (1-|\xi'|^2)(1+ |K|^2) - r^2}{1+|K|^2}.
\end{eqnarray}}
Observe that the support of $a_0$ is compactly contained in the interior of the set where $\tilde \rho = 1$.

We now quantize \eqref{decomposition} to see that 
\begin{eqnarray}
\label{operator factorization}
\frac{1}{1+K^2} h^2\tilde\Delta_\phi  = QJ + \t a_0(x',hD') - h\tilde e_1(x',hD') + h^2 \tilde e_0(x',hD')
\end{eqnarray}
where $\tilde e_1 = m_0 \t a_- \in S^{1}_1(\R^{n-1}), \tilde e_0 \in S^0_1(\R^{n-1})$, and $Q$ and $J$ are the operators with symbols $\xi_n - \tilde{a_-} + hm_0$ and $\xi_n - \tilde{a_+} + hm_0$ respectively. Observe that the $O(h)$ term in the composition formula for $QJ$ is killed by one of the $O(h)$ terms in \eqref{decomposition}.

Although this decomposition still gives us an $O(h)$ error, the symbol $\tilde e_1 $ vanishes when $|\xi'| =1$. In particular it vanishes on the characteristic set of $h^2\tilde\Delta_\phi$, and as the following lemma shows, it means that $h\tilde e_1(x',hD')\tilde{G}_{\phi}$ behaves one order of $h$ better than would be otherwise expected. This will help us with estimates later on. 

\begin{lemma} 
\label{composition estimate}
Let $\t E_1$ denote $\t e_1(x',hD')$. The operator $\tilde E_1 \tilde G_\phi$ is of the form
{\small \[\t E_1 \t G_\phi = (\t E_1 \t G_\phi)^c + Op_h(S^1_1 S^{-2}_1) + h\t e'_1(x',hD') \t G_\phi+ h Op_h(S^{-\infty}(\R^n)) \t G_\phi\]}
with $\t e'_1 \in S^0_1(\R^{n-1})$ and
\[(\t E_1 \t G_\phi)^c: L^2 \to_{h^0} H^k,\ \  (\t E_1 \t G_\phi)^c:L^{p'}\to_{h^{-1}} H^k \, \, \forall k\in \N.\]
%Furthermore, $\t E_1^* \t G_\phi^* : H^{-1} \to_{h^0} L^2$.
\end{lemma}
Here the notation $T:X \rightarrow_{h^m} Y$ indicates that the norm of the operator $T$ from $X$ to $Y$ is bounded by $O(h^m)$.

\begin{proof} We use the fact that $\t e_1$ takes value zero on the characteristic set of $\t G_\phi$. First write 
\[\tilde E_1 = \tilde e_1(x',hD')= Op_h (\t a_+^{-1}m_0 \t a_- \t a_+) = Op_h(\t a_+^{-1} m_0) Op_h(\t a_- \t a_+) + hOp_h \t e'_1(x',hD')\]
for some $\t e_1' \in S^{0}_1(\R^{n-1})$. Note that \[Op_h(\t a_-\t a_+) \t G_\phi = Op_h(\t a_-\t a_+) \gamma^*(\chi G_\phi) + Op_h(\t a_-\t a_+)\gamma^*((1-\chi)G_\phi)\] for some compactly supported smooth function $\chi(\xi)$ which is identically $1$ on the ball of radius $2$. %We note here that it suffices to prove the Lemma for $\gamma^*(\chi G_\phi)$ in place of $\t G_\phi$.
This means that 
\begin{eqnarray}
\label{pulled-back operators}
\t E_1 \t G_\phi = Op_h(\t a_+^{-1} m_0) Op_h(\t a_- \t a_+) \gamma^*(\chi G_\phi) + Op_h( S^1_1 S^{-2}_1 )+ h\t E'_1 \t G_\phi.
\end{eqnarray}
From Proposition \ref{green with no boundary} $\t G_\phi =\gamma^* G_\phi$ where $G_\phi$ is the Fourier multiplier $\frac{1}{\xi_n^2 + i2\xi_n + (1-|\xi'|^2)}$. We compute the $Op_h(\t a_- \t a_+) \gamma^*(\chi G_\phi)$ portion of this operator. 
\begin{eqnarray*}
Op_h(\t a_- \t a_+) \gamma^*(\chi(hD) G_\phi) &=& Op_h(\t a_- \t a_+) \gamma^*(\chi(hD) )\gamma^*(G_\phi)\\
&=&Op_h(\t a_- \t a_+) (\t \chi(x,hD) + hOp_h(S^{-\infty}))\gamma^*( G_\phi)
\end{eqnarray*}
where $\t \chi(x,\xi) \in S^{-\infty}$ is the pulled-back symbol of $\chi(\xi)$. Continuing by composing $Op_h(\t a_- \t a_+) \t \chi(x,hD)$ using symbol calculus,
 \begin{eqnarray*}
Op_h(\t a_- \t a_+) \gamma^*(\chi G_\phi) %&=&(Op_h(\t a_- \t a_+\t \chi) + hOp_h(S^{-\infty}))\gamma^*( G_\phi)\\
%&=&(\gamma^*Op_h( a_-  a_+ \chi) + hOp_h(S^{-\infty}))\gamma^*( G_\phi)\\
&=& \gamma^*(Op_h( a_-  a_+) (\chi(hD) G_\phi)) + hOp_h(S^{-\infty})\gamma^*( G_\phi)
\end{eqnarray*}
where $a_\pm(x,\xi) := (\gamma^* \t a_\pm)(x,\xi) = \t a_\pm(x', D\gamma^T \xi)$. We claim $Op_h( a_-  a_+) (\chi(hD) G_\phi))$ can be written as the sum of a $\Psi$DO with symbol in $S^{-\infty}(\R^n)$ and an operator 
{\Small\begin{eqnarray}
\label{characteristic part of product}
(Op_h( a_-  a_+) (\chi(hD) G_\phi))^c: L^2 \to_{h^0} H^k,\ \  (Op_h( a_-  a_+) (\chi(hD) G_\phi))^c:L^{p'}\to_{h^{-1}} H^k.\
\end{eqnarray}}
Inserting this into \eqref{pulled-back operators} would give us the Lemma.

We verify our claim. Observe that 

{\tiny\begin{eqnarray}\label{compose}
a_+a_- = \frac{(1-\rho_0)^2 (1-|\xi' + \xi_n K|^2)}{1 + |K|^2} + \frac{(1 +iK\cdot (\xi' + \xi_n K))^2 - (1-\rho_0)^2(1 + iK\cdot(\xi'+ \xi_n K))^2}{(1 + |K|^2)^2}.
\end{eqnarray}}
where $\rho_0(x',\xi) = \tilde\rho_0(\xi' + \xi_n K)$. Now $\tilde \rho_0(x',\xi') = 0$ if $|\xi'| \geq c'$ for some $c'<1$ and $K(x')$ is uniformly bounded. Therefore $\rho_0(x',\xi) = 0$ if 
\[\frac{1+c'}{2} \leq |\xi'| \leq 2-c'\ \ {\rm and}\ \ \ |\xi_n| \leq \frac{1-c'}{2 (\sup_{x'}|K(x')| +1)}.\]
Since the characteristic set of the the Fourier multiplier $\frac{1}{\xi_n^2 + i2\xi_n + (1-|\xi'|^2)}$ is compactly contained in this set, let $\chi_2(\xi)$ be a cutoff which is supported in this set and $1$ in a neighbourhood of the characteristic set and define
\[
(Op_h( a_-  a_+) (\chi(hD) G_\phi))^c := Op_h( a_-  a_+) (\chi(hD) \chi_2(hD)G_\phi).
\]
Write 
{\small\[(Op_h( a_-  a_+) (\chi(hD) G_\phi)) = (Op_h( a_-  a_+) (\chi(hD) G_\phi))^c + (Op_h( a_-  a_+) (\chi(hD) (1-\chi_2(hD))G_\phi).\]}
The second expression is $\Psi$DO of order $-\infty$ since it vanishes identically near the characteristic set and is therefore a compactly supported smooth multiplier. 

It remains to establish \eqref{characteristic part of product} for the part containing the characteristic set. Since $\rho_0$ vanishes identically on the support of $\chi_2$, it follows from \eqref{compose} that 
\[
Op_h(a_+a_-) \chi(hD)\chi_2(hD)G_\phi = Op_h\left(\frac{(1-|\xi' + \xi_n K|^2)}{(1+|K|^2)^2}\right)\chi(hD)\chi_2(hD) G_\phi.
\]

Note since $Op_h(\frac{(1-|\xi' + \xi_n K|^2)}{(1+|K|^2)^2})$ is a differential operator, proving \eqref{characteristic part of product} amounts to proving estimates for the operators $Op_h(\frac{(1-|\xi'|^2)\chi(\xi)}{{\xi_n^2 + i2\xi_n + (1-|\xi'|^2)}})$ and $Op_h(\frac{\xi_n\chi(\xi)\langle\xi'\rangle}{{\xi_n^2 + i2\xi_n + (1-|\xi'|^2)}})$. Crucially, these are both bounded Fourier multipliers with compact support and therefore map $L^2 \to H^k$ for all $k\in \N$ with norm $O(1)$. Therefore \[Op_h(a_+a_-) \chi_2(hD)\chi(hD) : L^2\to H^k\] with norm $O(1)$.

Moving on to the $L^{p'} \to H^k$ estimate we write $ \chi(hD)G_\phi =\chi(hD)G_\phi \chi_{100}(hD)$ where $\chi_{100}(\xi)$ is identically $1$ on the support of $\chi$. The estimate is then a result of the $L^2$ estimate and the fact that $\chi_{100}(hD): L^{p'} \to_{h^0} W^{1,p'} \hookrightarrow_{h^{-1}} L^2$ by Sobolev embedding.
\end{proof}
%To obtain the last estimate for the composition of the adjoints $\t E_1 ^* \t G_\phi^*$, observe that $\t G_\phi^* = \t G_{-\phi}$ has the same characteristic set and that $\t E_1^* = \overline{\t e_1}(x',hD') + hOp_h (S^0_1(\R^{n-1}))$. The estimate can be established by similar arguments as above.  

\end{subsection}

\end{section}

\begin{section}{Parametrices on the Half-Space}

In this section we construct parametrices for $h^2\tilde \Delta_{\phi}$ on the upper half space which give vanishing trace on the boundary. By a change of variables, we will later use these to build the Green's function of Theorem \ref{green's function}. Because the factoring in \eqref{operator factorization} contains a large error term $A_0$ at small frequencies, we will perform two separate constructions -- one for the large frequency case and one for the small frequency case. We split the two frequency cases by using the cutoff function $\t \rho:(\R^{n-1}) \rightarrow \R$ defined above equation \eqref{supp of rho tilde}.

\begin{subsection}{Parametrix for $h^2\tilde \Delta_\phi$ at large frequency}\label{large freq}
Let $\tilde G_\phi$ be the Green's function from Proposition \ref{green with no boundary}, and $J^+ := J^{-1} {\bf 1}_{\R^n_+}$ where $J^{-1}$ is defined as Proposition \ref{J inverse}. Let $\tilde\Omega\subset \R^n_+$ be a smooth bounded open subset of the upper half-space (with possibly a portion of the boundary intersecting $x_n = 0$) . We show that the operator 
\[
P_l :=(1- \tilde\rho (hD'))J^+ J \tilde G_\phi
\]
is a suitable parametrix for the operator $h^2\tilde\Delta_\phi$ in $\tilde \Omega$ at large frequencies.

We begin by showing that $P_l$ has mapping properties like those of $\t G_\phi$.
\begin{prop}
\label{estimates for P_l}
The map $P_l$ satisfies, for $\delta >0$, \[P_l : L^2_\delta(\R^n) \to_{h^{-1}} H^1_{\delta-1}(\R^n),\ \ P_l : L^{p'}(\R^n)\to_{h^{-2}} L^p(\R^n).\] Furthermore, $P_lv \in  H^1_{loc}(\R^n)$ with $P_lv \mid_{x_n = 0} = 0$ for all $v\in L^{p'}(\R^n)$.
\end{prop}
\begin{proof}
The weighted $L^2$ Sobolev norms come as a direct consequence of the mapping properties of $\t G_\phi$ and the fact that $J$, $J^{-1}$ arise from $S^k_0(\R^n)$.

For the mapping property from $L^{p'}(\R^n) \rightarrow L^p(\R^n)$, we split $\t G_\phi= \t G_\phi^c + (\t G_\phi - \t G_\phi^c)$ and observe
\[J^+J \t G_\phi^c : L^{p'} \xrightarrow[h^{-2}]{\t G_\phi^c} W^{k,p} \xrightarrow {J}W^{k-1,p}\xrightarrow{J^+} W^{1,p} \, \, \, \, \mbox{ and }\]
\[J^+J (\t G_\phi-\t G_\phi^c ): L^{p'} \xrightarrow{\t G_\phi -\t G_\phi^c} W^{2,p'} \xrightarrow {J}W^{1,p'} \hookrightarrow _{h^{-1}} L^2\xrightarrow{J^+} H^1 \hookrightarrow_{h^{-1}} L^p.\]

The above diagram also shows that $P_l v \in H^1_{loc}$ for all $v\in L^{p'}$ by omitting the last Sobolev embedding. The trace property then comes from the definition of $P_l$ and Proposition \ref{J inverse}.
% $\t G_\phi = \t G_\phi^c + (\t G_\phi - \t G_\phi^c)$. We already saw that the characteristic part maps into 
\end{proof}

We now have the following proposition for $P_l$. In the following statement we denote ${\bf 1}_{\t \Omega}$ to be the indicator function of $\t\Omega$. If $v\in L^r(\t\Omega)$ we use the notation ${\bf 1}_{\t\Omega} v$ to denote its trivial extension to a function in $L^r(\R^n)$.
\begin{prop} \label{P_l remainders}
Let $\t \Omega\subset \R^n_+$ be a bounded domain with $\partial \t \Omega \cap \{x_n = 0\} \neq \varnothing$. Denote by ${\bf 1}_{\t\Omega}$ the indicator function of $\t \Omega$. Then $P_l$ is a parametrix at large frequencies with vanishing trace on the boundary of the upper half space, in the sense that for all $v \in L^{p'}(\t \Om)$,
\[
{\bf 1}_{\t\Omega}h^2\t \Delta_\phi P_l {\bf 1}_{\t\Omega} v= (1 - \t \rho(hD') + R_l + hR_l')v,
\]
with
\[
P_l {\bf 1}_{\t \Omega} v\in H^1_{loc}(\R^n),\ \ P_l {\bf 1}_{\t \Omega} v\mid_{\partial\t\Omega \cap \{x_n = 0\}} = 0,
\]
where $R_l = {\bf 1}_{\t\Omega} R_l {\bf 1}_{\t\Omega}$ and $R_l' = {\bf 1}_{\t\Omega} R_l'{\bf 1}_{\t\Omega}$ have the estimates
\[
R_l : L^2 \to_{h} L^2,\ \ R_l : L^{p'}\to_{h^0} L^2,\ \ R_l' : L^r \to_{h^0} L^r,\ \ \ 1<r<\infty.
\]
\end{prop}

To prove this, we compute in the sense of distributions on $\R^n_+$ acting on $C^\infty_0(\R^n_+)$. Using \eqref{operator factorization}

{\Small\begin{eqnarray*}
h^2 \tilde \Delta_\phi P_l&=& (1-\tilde\rho) h^2\tilde \Delta_\phi J^+ J  \tilde G_\phi + [h^2\tilde\Delta_\phi, \tilde\rho] J^+J \tilde G_\phi \\
&=&(1-\tilde\rho)(1+K^2)(QJ + \t A_0 + h\tilde E_1 + h^2\tilde E_0) J^+J\tilde G_\phi + [h^2\tilde\Delta_\phi, \tilde\rho] J^+J \tilde G_\phi\\
&=&(1-\tilde\rho)(1+K^2)(Q {\bf 1}_{\R^n_+} J \tilde G_\phi) + (1-\tilde\rho)(1+K^2) \t A_0 J^+ J  \tilde G_\phi\\ &&+ h(1-\tilde\rho)(1+K^2)\tilde E_1 J^+J\tilde G_\phi + h^2(1-\tilde\rho)(1+K^2)\tilde E_0J^+ J \tilde G_\phi  + [h^2\tilde\Delta_\phi, \tilde\rho] J^+J  \tilde G_\phi
\end{eqnarray*}}
The first term requires some care. Testing this operator against $v\in C^\infty_0(\R^n)$ and $u \in C_0^\infty(\R^n_+)$ yields $\langle Q^*(1+K^2)(1-\tilde\rho)^* u,( {\bf 1}_{\R^n_+} J  \tilde G_\phi)v\rangle_{L^2(\R^n)}$. The operator $Q^*$ is a $\Psi$DO in the $\xi'$ direction but it is only a differential operator in the $\xi_n$ direction. Therefore the support does not spread in the $x_n$ direction. The operator $\tilde \rho(hD')$ is an operator only in the $\xi'$ direction and therefore does not spread support in the $x_n$ direction. As such $Q^*(1+K^2)(1-\tilde\rho)^* u$ vanishes in an open neighbourhood containing the closure of the lower half space and therefore for all $u\in C^\infty_0(\R^n_+)$ and $v \in C^\infty_0(\R^n)$,
 \[\langle Q^*(1+K^2) (1-\tilde\rho)^* u,( {\bf 1}_{\R^n_+} J \tilde G_\phi)v\rangle_{L^2(\R^n)} = \langle Q^* (1+K^2)(1-\tilde\rho)^* u,(J  \tilde G_\phi)v\rangle_{L^2(\R^n)}.\]
Therefore we may continue our computation:
{\Small\begin{eqnarray*}
h^2 \tilde \Delta_\phi P_l&=&(1-\tilde\rho)(1+K^2)(Q  J  \tilde G_\phi) + (1-\tilde\rho)(1+K^2) \t A_0 J^+ J \tilde G_\phi + h(1-\tilde\rho)(1+K^2)\tilde E_1 J^+J\tilde G_\phi \\&& + h^2(1-\tilde\rho)(1+K^2)\tilde E_0J^+ J \tilde G_\phi  + [h^2\tilde\Delta_\phi, \tilde\rho] J^+J \tilde G_\phi.\\
%&=&(1-\tilde\rho)(1+K^2)( Q  J  \tilde G_\phi) + (1-\tilde\rho) A_0 J^+ J   \tilde G_\phi + h(1-\tilde\rho)\tilde E_1 J^+J \tilde G_\phi \\&&+ h^2(1-\tilde\rho)\tilde E_0J^+ J  \tilde G_\phi  + [h^2\tilde\Delta_\phi, \tilde\rho] J^+J \tilde G_\phi \\
%&=& (1-\tilde\rho)(1+K^2)( Q  J  \tilde G_\phi) + (1-\tilde\rho) A_0 J^+ J   \tilde G_\phi + h(1-\tilde\rho)\tilde E_1 J^+J \tilde G_\phi \\&&+ h^2(1-\tilde\rho)\tilde E_0J^+ J   \tilde G_\phi  + [h^2\tilde\Delta_\phi, \tilde\rho] J^+J   \tilde G_\phi\\
\end{eqnarray*}}
At this juncture we invoke the factorization \eqref{operator factorization} again and plug the relation 
\[ h^2\tilde\Delta_\phi   -(1+K^2)(\t A_0 - h\tilde E_1 + h^2 \tilde E_0) =(1+K^2) QJ \]
into the first term. Since $h^2\tilde \Delta_{\phi} G_{\phi} = I$, we get 
%{\small\begin{eqnarray*}
%h^2 \tilde \Delta_\phi P_l&=&(1-\tilde\rho)( h^2\tilde\Delta_\phi   -(1+K^2)(\t A_0 - h\tilde E_1 + h^2 \tilde E_0)) \tilde G_\phi + (1-\tilde\rho)(1+K^2) \t A_0 J^+ J \tilde G_\phi \\&&+ h(1-\tilde\rho)(1+K^2)\tilde E_1 J^+J\tilde G_\phi  + h^2(1-\tilde\rho)(1+K^2)\tilde E_0J^+ J \tilde G_\phi  + [h^2\tilde\Delta_\phi, \tilde\rho] J^+J \tilde G_\phi\\
%&=&(1-\tilde\rho)( I   -(1+K^2)(\t A_0 - h\tilde E_1 + h^2 \tilde E_0) \tilde G_\phi) + (1-\tilde\rho)(1+K^2)\t A_0 J^+ J \tilde G_\phi \\&&+ h(1-\tilde\rho)(1+K^2)\tilde E_1 J^+J\tilde G_\phi  + h^2(1-\tilde\rho)(1+K^2)\tilde E_0J^+ J \tilde G_\phi  + [h^2\tilde\Delta_\phi, \tilde\rho] J^+J \tilde G_\phi\\
%\end{eqnarray*}}
for all $v\in C^\infty_0(\R^n)$ ,
\begin{eqnarray}
\label{large parametrix remainder gathered}
 h^2 \tilde \Delta_\phi P_l v = (1-\t \rho) v + R_1v+ R_2v +  [h^2\tilde\Delta_\phi, \tilde\rho] J^+J \tilde G_\phi v
 \end{eqnarray}
as a distribution on $\R^n_+$  (ie integrating against functions in $C^\infty_0(\R^n_+)$) where
\begin{eqnarray}
\label{R1}
R_1 = h(1-\tilde\rho) (1+K^2) \tilde E_1(1- J^+J) \tilde G_\phi \mbox{ and }
\end{eqnarray}
\begin{eqnarray}
\label{R2}
R_2 = (1-\tilde\rho)(1+K^2)(\t A_0 - \t A_0 J^+J + h^2\tilde E_0 - h^2\tilde E_0 J^+J)\tilde G_\phi.
\end{eqnarray}
In the following three lemmas, we claim that the remainder terms in \eqref{large parametrix remainder gathered} have the form of the remainders in Proposition \ref{P_l remainders}. The estimates for the terms in $R_2$ do not use the finer structures of $\tilde G_\phi$ while the estimates for terms in $R_1$ takes advantage of smallness of operators whose symbol is zero on the characteristic set of $\tilde G_\phi$.
\begin{lemma}
\label{commuting the laplacian}
\[ [h^2\t \Delta_\phi, \t \rho] J^+ J \t G_\phi = h^2 R_0'' + hR_0'\]
where $R_0': L^r\to L^r$, $R_0'' : L^2_\delta \to_{h^{-1}} L^2_{\delta-1}$, and $R_0'' : L^{p'}\to_{h^{-2}} L^{p}$.
\end{lemma}

\begin{lemma}
\label{R1 estimates}
The operator $R_1$ from \eqref{R1} can be written as $R_1 = R_1' + R_1''$ where
\[\|R_1'\|_{L^2 \to L^2} + h \|R_1'\|_{L^{p'} \to L^2} \leq Ch,\ \  \|R''_1\|_{L^2_{\delta}\to L^2_{\delta-1}} + h\|R_1''\|_{L^{p'}\to L^p} \leq Ch \]
\end{lemma}

\begin{lemma}
\label{R2 estimates}
The operator $R_2$ from equation \eqref{R2} maps $R_2 : L^2_\delta\to L^2_{-\delta}$ with norm $O(h)$ while $R_2 : L^{p'} \to L^p$ with norm $O(1)$.
\end{lemma}

The basic idea is that Lemma \ref{R2 estimates} follows from the smallness of $h^2 \t E_0$ and the fact that $A_0$ is supported only where $(1 - \t \rho)$ is zero, Lemma \ref{R1 estimates} follows from the good behaviour of $\t E_1$ given by Lemma \ref{composition estimate}, and Lemma \ref{commuting the laplacian} follows from the good behaviour of $\t G_{\phi}$ off of the support of $\t \rho$.

\begin{proof}[Proof of Lemma \ref{R2 estimates}] The terms involving $h^2 \t E_0$ can be estimated directly using the estimates for $\t G_\phi$ and $P_l$ in Propositions \ref{green with no boundary} and \ref{estimates for P_l}. The terms involving $\t A_0$ can be estimated by observing that since $\t \rho(\xi')$ is chosen to be identically $1$ in a neighbourhood of the support of $\t a_0(x',\xi')$, the operator \[(1-\t\rho(hD')(1+K^2) \t A_0 \in h^\infty Op_h(S^\infty(\R^{n-1})).\]
\end{proof}

\begin{proof}[Proof of Lemma \ref{R1 estimates}] We begin with the $h\tilde E_1 \tilde G_\phi$ term in \eqref{R1}. By Lemma \ref{composition estimate},
{\Small \begin{eqnarray}
\label{first term of R1}
h\t E_1 \t G_\phi = h(\t E_1 \t G_\phi)^c + hOp_h(S^1_1 S^{-2}_1) + h^2\t e'_1(x',hD') \t G_\phi+ h^2 Op_h(S^{-\infty}(\R^n)) \t G_\phi\end{eqnarray}}
with $\t e'_1 \in S^0_1(\R^{n-1})$ and
\[(\t E_1 \t G_\phi)^c: L^2 \to_{h^0} H^k,\ \  (\t E_1 \t G_\phi)^c:L^{p'}\to_{h^{-1}} H^k\ \  \forall k\in \N.\]
Our task is to sort the terms in this operator into the $R_1'$ bin and the $R_1''$ bin. The first term of \eqref{first term of R1} fits the mapping properties of objects in the $R'_1$ bin. The second term of \eqref{first term of R1} is a $\Psi$DO in the $S^1_1 S^{-2}_1$ class and therefore belongs to the $R_1''$ bin by Proposition \ref{sobolev mapping}. By Proposition \ref{green with no boundary}, the third term of \eqref{first term of R1} can be written as 
\[ h^2\t e'_1(x',hD') \t G_\phi =  h^2\t e'_1(x',hD') \t G_\phi^c +  h^2\t e'_1(x',hD')(\t G_\phi- \t G_\phi^c)\]
where $(\t G_\phi- \t G_\phi^c)$ is a $\Psi$DO with symbol in $S^{-2}_1(\R^n)$. The $\Psi$DO part is in the $R'_1$ bin since it behaves well on $W^{r,k}$ spaces and the estimate for $L^{p'}\to L^2$ can be obtained by doing semiclassical Sobolev embedding. For the characteristic part, $\t G_\phi^c$ takes $L^2_\delta \to_{h^{-1}} H^k_{\delta-1}$ and $L^{p'}\to_{h^{-2}} W^{k,p}$. Therefore the characteristic part belongs to the $R''_1$ bin. 

The reasoning for the third term of \eqref{first term of R1} also applies to the last term %\[h^2\mathrm{Op}_h(S^{-\infty}(\R^n)\t G_{\phi},\] 
of \eqref{first term of R1} and shows that it can also be sorted into the $R'_1$ and $R''_1$ bin.

We proceed next with the $h\t E_1 J^+ J \t G_\phi$ term of \eqref{R1}:
\begin{eqnarray}
\label{second term of R1}
h\t E_1 J^+ J \t G_\phi = hJ^+ J \t E_1\t G_\phi+h [\t E_1, J^{-1}] {\bf 1}_{\R^n_+} J\t G_\phi  + hJ^+ [J,\t E_1]\t G_\phi.\end{eqnarray}
In the above calculation we commuted $\t E_1$ and ${\bf 1}_{\R^n_+}$ since $\t E_1$ only acts in the $x'$ direction.

The first term above can be handled exactly the same as the $h\t E_1\t G_\phi$ term -- note that the argument for the terms in \eqref{first term of R1} shows that each of the constituent terms of $h\t E_1 \t G_{\phi}$ in \eqref{first term of R1} maps to $W^{1,r}$, and so applying ${\bf 1}_{\R^n_+}J$ presents no difficulty. For the first commutator term of \eqref{second term of R1}, Lemma \ref{composing with j^-1} and Proposition \ref{J inverse} show that $ [\t E_1, J^{-1}]  = hm(x,hD)$ for some \[m(x,hD): L^r\to L^r,\ \ m(x,hD) : L^2_\delta \to L^2_\delta.\]
Therefore, splitting $\t G_\phi$ into to its characteristic part $\t G_\phi^c$ and its $\Psi$DO part $\t G_\phi-\t G_\phi^c $ as in Proposition \ref{green with no boundary} we have
{\Small \[L^{p'} \xrightarrow{\t G_\phi -\t G_\phi^c} W^{2,p'} \xrightarrow{J} W^{1,p'} \hookrightarrow_{h^{-1}} L^2 \xrightarrow{{\bf 1}_{R^n_+}} L^2\xrightarrow [h]{[J^{-1},\t E_1]} L^2 ,\ \ \ L^{2} \xrightarrow{\t G_\phi -\t G_\phi^c} H^2\xrightarrow {J} H^1 \xrightarrow [h]{[\t E_1, J^{-1}] {\bf 1}_{\R^n_+}} L^2\]}
and so $h [\t E_1, J^{-1}] {\bf 1}_{\R^n_+} J(\t G_\phi- \t G_\phi^c)$ belongs to $R'_1$ bin.
For the characteristic part
{\Small \[L^{p'} \xrightarrow[h^{-2}]{\t G_\phi^c} W^{k,p} \xrightarrow{J} W^{k-1,p}  \xrightarrow{{\bf 1}_{R^n_+}}L^p\xrightarrow [h]{[J^{-1},\t E_1]} L^p ,\ \ \ L_\delta^{2} \xrightarrow{\t G_\phi^c} H_{\delta-1}^k\xrightarrow {J} H^{k-1}_{\delta-1} \xrightarrow [h]{[\t E_1, J^{-1}] {\bf 1}_{\R^n_+}} L_{\delta-1}^2\]}
and therefore $h [\t E_1, J^{-1}] {\bf 1}_{\R^n_+} J\t G_\phi ^c$ belongs to the $R''_1$ bin.

For the $[J,\t E_1]\t G_\phi$ term, splitting $\t G_\phi$ into to its characteristic part $\t G_\phi^c$ and its $\Psi$DO part $\t G_\phi-\t G_\phi^c $ we have 
{\small \[L^{p'} \xrightarrow{\t G_\phi -\t G_\phi^c} W^{2,p'} \xrightarrow [h]{[J,\t E_1]} W^{1,p'} \hookrightarrow _{h^{-1}} L^2 \xrightarrow {J^+} H^1,\ \ \ L^{2} \xrightarrow{\t G_\phi -\t G_\phi^c} H^2\xrightarrow [h]{[J,\t E_1]} H^1 \xrightarrow {J^+} H^1.\]}
Therefore $hJ^+ [J,\t E_1] (\t G_\phi - \t G^c)$ belongs to the $R_1'$ bin. For the part with characteristic set,  $J^+ [J,\t E_1] \t G^c$ behaves like
{\small \[L^{p'} \xrightarrow[h^{-2}]{\t G_\phi^c} W^{k,p} \xrightarrow [h]{[J,\t E_1]} W^{k-1,p} \xrightarrow {J^+} W^{1,p},\ \ \ L^{2}_\delta \xrightarrow[h^{-1}]{\t G_\phi^c} H_{\delta-1}^k\xrightarrow [h]{[J,\t E_1]} H_{\delta-1}^{k-1} \xrightarrow {J^+} H_{\delta-1}^1\]}
and therefore $hJ^+ [J,\t E_1] \t G^c$ belongs to $R''_1$ bin.
\end{proof}

% In the following computation, $ a(x',hD) = O(h)$ means \[a(x',hD) = h m(x',hD),\ \ \|m(x',hD)\|_{L^r\to L^r} \leq C\ \ \ \|m(x',hD)\|_{L^2_\delta\to L^2_\delta} \leq C.\]

%By Lemma \ref{composing with j^-1} and Theorem \ref{J inverse}, $\tilde E_1 J^{-1} = \tilde E_1 j^{-1}(x',hD) + O(h)$. By Lemma \ref{composing with j^-1} we may commute $\tilde E_1$ and $j^{-1}(x',hD)$ to obtain $\tilde E_1 J^{-1} =  j^{-1}(x',hD)  \tilde E_1+ O(h)$.  We are allowed to commute $\tilde E_1$ and ${\bf 1}_{\R^n_+}$ by the fact that $\tilde E_1$ is only a $\Psi$DO in the $x' \in \R^{n-1}$ variable and ${\bf 1}_{R^n_+}$ is constant along $x'$. Therefore, $\tilde E_1 J^+ = j^{-1} (x', hD){\bf 1}_{\R^n_+} \tilde E_1 + O(h)$ and this gives 
 %\[\tilde E_1 J^+J \tilde G_\phi =  j^{-1} (x', hD){\bf 1}_{\R^n_+} \tilde E_1 J\tilde G_\phi + O(h)J\tilde G_\phi.\]
 %Lastly we would like to commute $\tilde E_1$ and $J$. For this we use the explicit expression
 %\[ J = hD_n + F(x',hD') + hF_0(x',hD')\] to obtain $[\tilde E_1, J]= [\tilde E_1, F(x',hD')] = O(h)$. \qed

%Finally we treat the commutator term with Laplacian in \eqref{large parametrix remainder gathered}. 

\begin{proof}[Proof of Lemma \ref{commuting the laplacian}] We have
\[[h^2\t \Delta_\phi, \t \rho] J^+ J \t G_\phi= [K^2, \t \rho] h^2D_n^2 J^+J\t G_\phi - 2[K\cdot hD_{x'}, \t \rho] hD_n J^+J\t G_\phi.\]
Some care will be needed in treating the term involving $h^2D_n^2$ hitting $J^+ = J^{-1} {\bf 1}_{\R^n_+}$. We are only considering the expressions as maps to distributions on $\R^n_+$, so for all $u\in C^\infty_0(\R^n_+)$ and $v \in  C^\infty_0(\R^n)$,
\[ 
\langle hD_nu, hD_n J^{-1}{\bf 1}_{\R^n_+}v\rangle = \langle hD_n u, (1- FJ^+) v\rangle= \langle u, hD_n v - F v - FJ^+ v\rangle.
\]
Here we used the fact that $J = hD_n + F(x',hD')$ for some $F(x',\xi') \in S^{1}_1(\R^{n-1})$ and the tangential operator $F(x',hD')$ commutes with the indicator function of the upper half-space.

Combining the two expressions we obtain
{\Small\begin{eqnarray}
\label{commutator of laplace simplified}
[h^2\t \Delta_\phi, \t \rho] J^+ J \t G_\phi = [K^2,\t \rho] (hD_n - F - FJ^+) J\t G_\phi - 2[K\cdot hD_{x'}, \t\rho] (1-FJ^+)J\t G_\phi.
\end{eqnarray}}
 We decompose $\t G_\phi$ in \eqref{commutator of laplace simplified} into its $\Psi$DO part and its characteristic part as stated in Proposition \ref{green with no boundary}. The part of \eqref{commutator of laplace simplified} containing the $\Psi$DO is a bounded map from $L^r\to L^r$ with a gain in $h$ obtained from the commutator. Therefore, the part containing the $\Psi$DO belongs to the $hR'_0$ bin.  

For the part containing the characteristic set, we expand $[h^2\t \Delta_\phi, \t \rho] J^+ J \t G_\phi^c$ as
%{\Small
\begin{eqnarray*}
& &  [K^2,\t \rho] \t\rho_1 (hD_n - F - FJ^+) J\t G_\phi^c - 2[K\cdot hD_{x'}, \t\rho] \t\rho_1(1-FJ^+)J\t G_\phi^c\\ 
&+&  [K^2,\t \rho] (1-\t\rho_1) (hD_n - F - FJ^+) J\t G_\phi^c - 2[K\cdot hD_{x'}, \t\rho] (1-\t\rho_1) (1-FJ^+)J\t G_\phi^c
\end{eqnarray*}%}
where $\t\rho_1(\xi')$ is chosen to be identically $1$ in a neighbourhood compactly containing the support of $\t \rho$ but vanishes identically in a neighbourhood of $|\xi'| = 1$. By disjoint support, $ [K^2,\t \rho] (1-\t\rho_1)$ and $[K\cdot hD_{x'}, \t\rho] (1-\t\rho_1) $ both belong to $h^\infty S^{-\infty} (\R^{n-1})$. Since $\t G_\phi^c : L^2_{\delta} \to_{h^{-1}} H^k_\delta$ and $L^{p'} \to _{h^{-2}} W^{k,p}$, the second line in the above expression for $[h^2\t \Delta_\phi, \t \rho] J^+ J \t G_\phi^c $ can be sorted into the $h^2 R''_0$ bin.

The only thing remaining is to treat the terms on the support of $\t\rho_1$. We will treat the first term and the second term is dealt with in the same manner. We commute $\t\rho_1(hD')$ so that it appears next to $\t G_\phi^c$:
\[[K^2,\t \rho] \t\rho_1 (hD_n - F - FJ^+) J\t G_\phi^c = [K^2,\t \rho]  (hD_n - F - FJ^+) J \t \rho_1\t G_\phi^c + h^2 R''_0.\]
We are able to throw all the commutator terms with $\t \rho_1$ into the $h^2 R''_0$ bin by using Proposition \ref{J inverse}, Lemma \ref{composing with j^-1}, and Proposition \ref{sobolev composition} in conjunction with the mapping properties of $\t G_\phi^c$ given by Proposition \ref{green with no boundary}. Since $\t \rho_1(\xi')$ vanishes identically near $|\xi'| = 1$, Lemma \ref{disjoint from char} asserts that,  \[\t \rho_1\t G_\phi^c = Op_h(S^{-\infty}(\R^n)) + hm(x',hD) \t G_\phi^c\] for some $m(x,\xi)\in S^{-\infty}(\R^n)$ and 
therefore \[[K^2,\t \rho] \t\rho_1 (hD_n - F - FJ^+) J\t G_\phi^c = hR_0' + h^2 R_0''.\]
\end{proof}

\begin{proof}[Proof of Proposition \ref{P_l remainders}]
The estimates for $R_l$ and $R_l'$ come from Lemmas \ref{commuting the laplacian}, \ref{R1 estimates}, and \ref{R2 estimates} in conjunction with \eqref{large parametrix remainder gathered}. The trace property of the operator $P_l {\bf 1}_{\t \Omega}$  on $\partial\t\Omega \cap \{x_n = 0\}$ is a result of Proposition \ref{estimates for P_l}. Note that the $L^2$ bounds in Proposition \ref{P_l remainders} are unweighted because of the conjugation with indicator functions of $\tilde{\Om}$. 
\end{proof}

%\begin{cor}
%\label{P_l remainders}
%Let $\t \Omega\subset \R^n_+$ be a bounded set with possibly overlapping boundary components with $\{x_n = 0\}$. Denote by ${\bf 1}_{\t\Omega}$ the indicator function then $P_l$ is a parametrix in the sense that 
%\[ {\bf 1}_{\t\Omega} h^2\t \Delta_\phi P_l {\bf 1}_{\t\Omega}= {\bf 1}_{\t\Omega} +  R_l + hR_l',\ \ P_l\mid_{\partial\t\Omega \cap \{x_n = 0\}} = 0\]
%where $R_l = {\bf 1}_{\t\Omega} R_l {\bf 1}_{\t\Omega}$ and $R_l' = {\bf 1}_{\t\Omega} R_l'{\bf 1}_{\t\Omega}$ have the estimates
%\[R_l : L^2 \to_{h} L^2,\ \ R_l : L^{p'}\to_{h^0} L^2,\ \ R_l' : L^r \to_{h^0} L^r,\ \ \ 1<r<\infty\]
%\end{cor}

\end{subsection}
\begin{subsection}{Parametrix for $h^2\t\Delta_\phi$ at Small Frequency}\label{small freq}

Here we want to look for a parametrix for $h^2 \tilde{\Delta}_{\ph}$ at low frequencies. We begin by defining $p(x',\xi)$ to be the symbol of $h^2\tilde{\Delta}_{\ph}$:
\[
p(x',\xi) := (1+K^2)\xi_n^2 -2\xi_n(i-\xi'\cdot K) -(1- |\xi'|^2).
\]
Now  define 
\[
P_s := \frac{\t \rho}{p}(x',hD).
\]
The following proposition says that $P_s$ inverts $h^2\tilde{\Delta}_{\phi}$ at small frequencies, up to an $O(h)$ error.

\begin{prop}
\label{small freq parametrix} $P_s$ is a bounded operator $P_s : L^r \to W^{2,r}$ for all $r\in (1,\infty)$. Moreover for all $r\in (1,\infty)$.
\[
h^2 \t \Delta_\phi P_s = \t\rho + hR_s
\]
for some $R_s : L^r\to L^r$ bounded uniformly in $h$.
\end{prop}

\begin{proof}
We want to use the symbol calculus developed in Section 2. However, we have the complication that 
%\lt{just realized that $p$ is both the polynomial and the index of $L^p$ space. Maybe we write out explicitly $p(x',\xi)$ whenever we talk about the polynomial.}
$1/p(x',\xi)$ is not a proper symbol, because of the zeros of $p(x',\xi)$. Therefore it is not immediately evident that $\t \rho/p(x',\xi)$ lies in the symbol class $S^{-\infty}S^{-2}_1$, as we would want.

We can remedy this by writing
\[
\t \rho(\xi')/p (x',\xi)  = (1 - \chi_{100}(\xi))\t \rho(\xi')/p(x',\xi) + \chi_{100}(\xi)\t \rho(\xi')/p(x',\xi)
\]
where $\chi_{100}(\xi) \in S^{-\infty}(\R^n)$ is a smooth cutoff function supported only for $|\xi| < 100$, and identically one in the ball $|\xi| \leq 50$.  

Now note that by \eqref{supp of rho tilde}, $p(x',\xi)$ is properly elliptic on the support of $\t \rho(\xi')$, so $\chi_{100}(\xi) \t \rho(\xi')/p (x',\xi) \in S^{-\infty}(\R^n)$. Moreover, since the characteristic set of $p(x',\xi)$ lies well inside the set where $\chi_{100} \equiv 1$, we have that $(1 - \chi_{100}(\xi))/p(x',\xi) \in S^{-2}_1(\R^n)$.  

Therefore $P_s$ can be understood as the sum of two operators, one of which is in the symbol class $S^{-\infty}(\R^n)$ and the other of which is in the symbol class $S^{-\infty}S^{-2}_1(\R^n)$. Then Proposition \ref{sobolev mapping} asserts that $P_s : L^r \to W^{2,r}$ is a bounded operator and Proposition \ref{sobolev composition} asserts that  
\[
h^2\tilde{\Delta}_{\ph}\mathrm{Op}_h\left( \frac{\t \rho}{p} \right) = \mathrm{Op}_h((1 - \chi_{100})\t \rho) + \mathrm{Op}_h(\chi_{100}\t \rho) + hR_{-1} = \mathrm{Op}_h (\t \rho) + hR_s 
\]
as we wanted.  

\end{proof}

It turns out that our small frequency parametrix preserves support in the $x_n$ direction.

\begin{prop}\label{SmallFrequencySupport}
Suppose $v \in L^r(\R^n)$, with $1 < r < \infty$, and ${\supp}(v)$ is contained in the closure of $\R^n_+$. Then both ${\supp} (P_sv)$ and ${\supp} (R_sv)$ are contained in $\bar\R^n_+$, where $R_s$ is the operator from Proposition \ref{small freq parametrix}. In particular, $P_s v\mid_{x_n = 0} = 0$ if ${\supp} (v) \subset \bar\R^n_+$.
\end{prop}

\begin{proof}
Let $v \in C^{\infty}_c(\R^n)$.  Then 
\begin{equation}\label{Oprhopph}
\mathrm{Op}_h\left( \frac{\t \rho}{p} \right)v(x) = h^{-n}\int_{\R^n} \frac{\t \rho(\xi')}{p(x',\xi)}\hat{v}(\xi)e^{i\xi\cdot x/h} \, d\xi.
\end{equation}
We split the integral on the right into $x'$ and $x_n$ variables and get
\[
h^{-n}\int_{\R^{n-1}}e^{i\xi' \cdot x'/h} \int_{-\infty}^{\infty}\frac{\t \rho(\xi')}{p(x',\xi)}\hat{v}(\xi)e^{i\xi_n x_n/h} \, d\xi_n \, d\xi'.
\]
Consider the inner integral 
\[
\int_{-\infty}^{\infty}\frac{\t \rho(\xi')}{p(x',\xi)}\hat{v}(\xi)e^{i\xi_n x_n/h} \, d\xi_n.
\]
For fixed $\xi'$ and $x'$, we can write the Fourier transform of $v$ in the $\xi_n$ variable explicitly to get
%\[
%\int_{-\infty}^{\infty}\frac{\t \rho(\xi')}{p(x',\xi)}\int_{-\infty}^{\infty}\mathcal{F}_{x'}f(\xi',s)e^{-is\xi_n/h}\, ds \, e^{i\xi_n x_n/h} \, d\xi_n,
%\]
%which can be rewritten as 
\begin{equation}\label{InnerIntegral}
\int_{-\infty}^{\infty}\int_{-\infty}^{\infty}\frac{\t \rho(\xi')\mathcal{F}_{x'}v(\xi',s)e^{i\xi_n (x_n-s)/h}}{p(x',\xi)} \, d\xi_n \, ds.
\end{equation}
We want to evaluate the inner integral using the residue calculus. 
%\begin{equation}\label{ResidueIntegral}
%\int_{-\infty}^{\infty}\frac{\t \rho(\xi')\mathcal{F}_{x'}f(\xi',s)e^{i\xi_n (x_n-s)/h}}{p(x',\xi)} \, d\xi_n.
%\end{equation}
%To evaluate this, we will need to use the residue calculus.  
Since $e^{i\xi_n (x_n-s)/h}$ is analytic, we need to understand the zeros of $p(x',\xi)$ as a polynomial in $\xi_n$. Factoring, we have
\[
p(x',\xi) = -(1 + |K|^2)(\xi_n - a_+)(\xi_n - a_-)
\]
where
\[
a_{\pm} = i\frac{1+iK\cdot \xi' \pm \sqrt{(1+iK\cdot \xi')^2 - (1-|\xi'|^2)(1 +|K|^2)}}{1+|K|^2}.
\]
Therefore $p(x',\xi)$, viewed as a polynomial in $\xi_n$, has two roots: $a_+$ and $a_-$. Since we are taking the standard branch of the square root, it follows that $a_+$ has positive imaginary part.  Meanwhile, if the imaginary part of $a_-$ vanishes, then by proper choice of $\xi_n$, the factor $(\xi_n - a_-)$ can be made to vanish. On the other hand $\t \rho(\xi')$ is defined to have support only where $p$ is elliptic, and so the imaginary part of $a_-$ does not approach zero on the support of $\t \rho(\xi')$. Moreover $a_-$ has positive imaginary part when $\xi' = 0$, and it is continuous in $x'$ and $\xi'$ except when  $i\sqrt{(1+iK\cdot \xi')^2 - (1-|\xi'|^2)(1 +|K|^2)}$ is entirely real, and so $a_{-}$ also lies in the upper half of the complex plane for all $x'$ and $\xi'$.  

%Now if we examine the integrand in \eqref{InnerIntegral} we see that when $(x_n - s)\mathrm{Im}\xi_n \geq 0$, then the integrand goes to zero as $|\xi_n| \rightarrow \infty$. Therefore we can evaluate the $d\xi_n$ integral in \eqref{InnerIntegral} as a contour integral over the upper half complex plane when $x_n -s \geq 0$, and as a contour integral over the lower half complex plane when $x_n-s < 0$.  

%Since the integrand has no poles in the lower half complex plane, we find that the inner integral in \eqref{InnerIntegral} equals zero when $x_n - s < 0$. On the other hand, when $x_n -s \geq 0$, we pick up the residues from both poles $a_+(x',\xi')$ and $a_-(x',\xi')$.  Therefore the inner integral in \eqref{InnerIntegral} equals 
%\[
%2\pi i \frac{\t \rho(\xi')\mathcal{F}_{x'}f(\xi',s)e^{ia_- (x_n-s)/h}}{(1+|K|^2)(a_+ - a_-)}  - 2\pi i \frac{\t \rho(\xi')\mathcal{F}_{x'}f(\xi',s)e^{ia_+ (x_n-s)/h}}{(1+|K|^2)(a_+ - a_-)},
%\]
%for $x_n - s \geq 0$, at least when $a_+ \neq a_-$.  
%Then \eqref{InnerIntegral} becomes
Therefore evaluating the inner integral of \eqref{InnerIntegral} using the residue calculus over the appropriate contours, we get  
\[
2\pi i \int_{-\infty}^{x_n}\frac{\t \rho(\xi')\mathcal{F}_{x'}v(\xi',s)(e^{ia_- (x_n-s)/h}-e^{ia_+ (x_n-s)/h})}{(1+|K|^2)(a_+ - a_-)} \, ds
\]
at least for $a_+ \neq a_-$. Note that since $a_{\pm}$ both have positive imaginary part on the support of $\t \rho(\xi')$, this integral converges. Now
{\Small\begin{equation}\label{SupportForm}
\mathrm{Op}_h\left( \frac{\t \rho}{p} \right)f(x) = 2\pi i h^{-n}\int_{\R^{n-1}}e^{i\xi' \cdot x'/h} \int_{-\infty}^{x_n}\frac{\t \rho(\xi')\mathcal{F}_{x'}v(\xi',s)(e^{ia_- (x_n-s)/h}-e^{ia_+ (x_n-s)/h})}{(1+|K|^2)(a_+ - a_-)} \, ds \, d\xi'.
\end{equation}}
At first glance this integral may have issues with convergence when $a_+ -a_- \rightarrow 0$. However, on the set where $a_+ = a_-$, the residue calculus tells us that the integral vanishes, and near this set we have
%\lt{you had multiplication. I corrected it to be division.}
\[
\lim_{a_+ - a_- \rightarrow 0}\frac{e^{ia_- (x_n-s)/h}-e^{ia_+ (x_n-s)/h}}{(1+|K|^2)(a_+ - a_-)} = \frac{i(x_n-s)}{h(1 + |K|^2)}e^{ia_- (x_n-s)/h}.
\]
Therefore the integral on the right side of \eqref{SupportForm} converges, and so this provides an honest representation of $P_s = \mathrm{Op}_h(\t \rho/p)$, at least when $v \in C^{\infty}_{c}(\R^n)$.  Note that we are not claiming that this integral proves $L^r$ boundedness: the non-smoothness of $a_{\pm}$ makes this non-obvious. Rather, we want to use this representation of the operator to prove the support property. If $v \in C^{\infty}_0(\R^n)$ is supported only in the upper half space $x_n > 0$, it is clear from \eqref{SupportForm} that 
\begin{equation}\label{SupportProp}
P_s v(x',x_n) = 0 \mbox{ for } x_n \leq 0.
\end{equation}
Now from Proposition \ref{small freq parametrix} we have
\[
\|P_s v\|_{W^{2,r}(\R^n)} \leq \|v\|_{L^r(\R^n)}, 
\]
and it follows from the trace theorem that for any fixed $x_n$,
\begin{equation}\label{TraceProperty}
h\|P_s v(\cdot, x_n)\|_{W^{1,r}(\R^{n-1})} \leq \|v\|_{L^r(\R^n)}.
\end{equation}
Therefore if $v \in L^r(\R^n)$ is supported only in the upper half space, we can approximate it with $C^{\infty}_0$ functions supported in the upper half space and use the support property for those functions, together with \eqref{TraceProperty}, to conclude that 
\[
\mathrm{Op}_h\left( \frac{\rho}{p} \right)v(x',x_n) = 0
\]
for $x_n \leq 0$. This shows that $P_s$ has the desired support property. The support property for $R_s$ then follows from writing 
\[
h^2 \t \Delta_\phi P_s - \t \rho(h D') =  hR_s
\]
and noting that every operator on the left hand side of this equation has the desired support property.

\end{proof}

\end{subsection}

\end{section}

\begin{section}{Dirichlet Green's function and Carleman estimates}

\begin{subsection}{Green's Function For Single Graph Domains}\label{graph case}

%{\bf I screwed up below and used the notation $G_\phi$ for the green's function with boundary condition when earlier I used the same symbol to denote the constant coefficient multiplier used by Sylvester-Uhlmann. The notation will need to be changed at some point.}
By combining Sections \ref{large freq} and \ref{small freq} %in conjunction with the estimates in Corollary \ref{P_l remainders} 
we see that ${\bf 1}_{\t \Omega} (P_s + P_l) {\bf 1}_{\t \Omega}$ is a parametrix for the operator $h^2\t\Delta_\phi$ in the domain $\t \Omega$. As one expects, this parametrix can be modified into a Green's function. 

In this section we consider domains with a component of the boundary which coincides with the graph of a function. In particular, let $\Omega$ be a bounded domain in $\R^n$, and suppose $f \in C_0^{\infty}(\R^{n-1})$ such that $\Omega$ lies in the set $\{x_n > f(x')\}$, with a portion of the boundary $\Gamma \subset \partial \Omega$ lying on the graph $\{x_n = f(x')\}$. Denote by $\gamma$ the change of variable $(x',x_n) \mapsto (x', x_n - f(x'))$.

\begin{prop}
\label{graph green}
There exists a Green's function $G_\Gamma$ which satisfies the relation $\langle h^2\Delta_\phi^* u, G_\Gamma f\rangle = \langle u, f\rangle$ for all $u\in C^\infty_0(\Omega)$ and is of the form 
\[\gamma_*G_\Gamma = {\bf 1}_{\t \Omega}(P_s + P_l) {\bf 1}_{\t \Omega}  (I + R) \] with $R$ obeying the estimates
\[ R : L^{p'} (\t\Omega) \to_{h^0} L^2(\t \Omega),\ \ \ R : L^{2} (\t\Omega) \to_{h} L^2(\t \Omega).\]
The Green's function $G_\Gamma$ satisfies the estimates
\[G_\Gamma: L^2(\Omega)\to_{h^{-1}} L^2(\Omega),\ G_\Gamma : L^{p'}(\Omega)\to_{h^{-2}} L^p(\Omega).\] Furthermore, $G_\Gamma v\in H^1(\Omega)$ for all $v \in L^{p'}$  and $G_\Gamma v \mid_{\Gamma} = 0$.
\end{prop}
\begin{proof}
Change coordinates $(x',x_n) \mapsto (x', x_n - f(x'))$ so that $\t \Gamma\subset \{x_n = 0\}$ and let $\tilde \Delta_\phi$ be the pulled-back conjugated Laplacian. All equalities below are in the sense of distributions in $\t \Omega$. By Proposition \ref{P_l remainders} and Proposition \ref{small freq parametrix}, for any $v\in L^{p'}(\t \Omega)$,
\[\langle h^2\tilde\Delta_\phi^* u, {\bf 1}_{\t \Omega} (P_s + P_l) {\bf 1}_{\t \Omega} v \rangle = \langle u, v + (hR_s + hR_l' + R_l ) v\rangle\ \ \ \forall u\in C^\infty_0(\t\Omega)\]
with $R_s$ and $R_l'$ mapping $L^r \to L^r$ with no loss in $h$ while
\[R_l : L^2 \to_{h} L^2,\ \ R_l : L^{p'}\to_{h^0} L^2.\]
Let $S :L^r\to L^r$ denote the inverse of $(1 + hR_l' + h R_s)$ by Neumann series. Then in $\t \Omega$ we have
\[h^2\tilde\Delta_\phi {\bf 1}_{\t \Omega}(P_s + P_l){\bf 1}_{\t \Omega}S = I + R_l S,\]
with $R_lS : L^2\to_{h} L^2$ while $R_lS : L^{p'} \to_{h^0} L^2$. Therefore, for all $v \in L^{p'}(\t \Omega)$ the Neumann series
\[(1 + R_l S)^{-1}v := v- \sum\limits_{k=0}^\infty (-R_lS)^k (R_lS)v \in L^{p'}\]
is well-defined and the series converge in $L^2(\t\Omega)$.
The operator ${\bf 1}_{\t \Omega}(P_s + P_l) {\bf 1}_{\t \Omega} S(1 + R_l S)^{-1} $ is then a right inverse of $h^2\t \Delta_\phi$ in $\t \Omega$. By defining 
\[G_\Gamma := \gamma^*{\bf 1}_{\t \Omega}(P_s + P_l) {\bf 1}_{\t \Omega} S(1 + R_l S)^{-1}\]
one obtains the Green's function in the original coordinates.

For the estimates on $G_\Gamma$ and for verifying the trace it is more convenient to work with the operator ${\bf 1}_{\t \Omega}(P_s + P_l) {\bf 1}_{\t \Omega} S(1 + R_l S)^{-1}$ and deduce the analogous properties for $G_\Gamma$. We first check that ${\bf 1}_{\t \Omega}(P_s + P_l) {\bf 1}_{\t \Omega} S(1 + R_l S)^{-1}  v\in H^1(\t\Omega)$ for all $v\in L^{p'}$ and that the trace vanishes on $\t \Gamma \subset\{x_n = 0\}$. By Proposition \ref{estimates for P_l} the operator $P_l$ maps $L^{p'}$ into $H^1_{loc}$ has vanishing trace on $\{x_n = 0\}$. By Proposition \ref{small freq parametrix} $P_s v$ is an element of $W^{2,p'} (\R^n) \hookrightarrow H^1(\R^n)$ which vanishes in $\{x_n <0\}$ if $v\in L^{p'}(\R^n)$ vanishes in $\{x_n <0\}$. Therefore we conclude that ${\bf 1}_{\t \Omega}(P_s + P_l) {\bf 1}_{\t \Omega} S(1 + R_l S)^{-1}  v\in H^1(\t\Omega)$ has trace zero on $\t \Gamma$ for all $v \in L^{p'}(\t \Omega)$ and thus $G_\Gamma$ has vanishing trace on $\Gamma$.

To verify the mapping properties of ${\bf 1}_{\t \Omega}(P_s + P_l) {\bf 1}_{\t \Omega} S(1 + R_l S)^{-1} $ write
{\small\begin{eqnarray*}{\bf 1}_{\t \Omega}(P_s + P_l) {\bf 1}_{\t \Omega} S(1 + R_l S)^{-1}  &=&{\bf 1}_{\t \Omega} (P_s + P_l) {\bf 1}_{\t \Omega} S(I- \sum\limits_{k=0}^\infty (R_lS)^k (R_lS) )\\
&=&{\bf 1}_{\t \Omega} (P_s + P_l) {\bf 1}_{\t \Omega} S - {\bf 1}_{\t \Omega} (P_s + P_l) {\bf 1}_{\t \Omega} S\sum\limits_{k=0}^\infty (R_lS)^k (R_lS)
\end{eqnarray*}}
Since $S: L^r \to L^r$, inserting an $L^2(\t \Omega)$ function would yield, by Propositions \ref{estimates for P_l} and \ref{small freq parametrix},
 an $H^1$ function with a loss of $h^{-1}$ in the first term and no loss in the second. For mappings from $L^{p'}$ we only need to concern ourselves with the first term since the Neumann sum maps $L^{p'}\to L^2$ with no loss in $h$ and we can refer to the $L^2$ estimate for ${\bf 1}_{\t \Omega} (P_s + P_l) {\bf 1}_{\t \Omega} S$.

The mapping properties are then verified by observing, due to Proposition \ref{estimates for P_l},
\[{\bf 1}_{\t \Omega} P_l {\bf 1}_{\t \Omega} S: L^{p'} \xrightarrow{S} L^{p'}\xrightarrow{{\bf 1}_{\t\Omega}}L^{p'}\xrightarrow[h^{-2}]{P_l} L^p\xrightarrow{{\bf 1}_{\t\Omega}}L^{p}.\]
And due to Proposition \ref{small freq parametrix},
\[{\bf 1}_{\t \Omega} P_s {\bf 1}_{\t \Omega} S: L^{p'} \xrightarrow{S} L^{p'}\xrightarrow{{\bf 1}_{\t\Omega}}L^{p'}\xrightarrow{P_s} W^{2,p'} \hookrightarrow_{h^{-2}} L^p\xrightarrow{{\bf 1}_{\t\Omega}}L^{p}.\]
%For the trace property, since $G_\Gamma = \gamma^*{\bf 1}_{\t \Omega}(P_s + P_l) {\bf 1}_{\t \Omega} S(1 + R_l S)^{-1}$ we can conclude that $G_\Gamma$ has vanishing trace on $\Gamma\subset \partial \Omega$.
\end{proof}

This finishes the proof of Theorem \ref{green's function} in the case when $\Gamma$ lies in a single graph. In the next section we move on to the general case.

\end{subsection}
\begin{subsection}{Proof of Theorem \ref{green's function} - Dirichlet Green's Function}\label{general case}
%In this section we prove Theorem \ref{green's function}. % Let $\Omega \subset \R^n$ be a bounded open subset with smooth boundary. Let $\omega \in \R^n$ be a unit vector and $\Gamma \subset \partial\Omega$ be an open set compactly contained in $\{ x\in \partial\Omega \mid \nu \cdot \omega >0\} $. Then by a rigid rotation one can write each connected piece $\Gamma_j$ as a portion of the graph of a compactly supported smooth function $x_n = f_j(x')$.
%\begin{prop}
%\label{green general domain}
%Let $\Omega$ and $\Gamma_j$ be the sets described above, with $\Gamma_j$ lying on the graph $\{x_n = f_j(x')\}$. Let $\Delta_\phi := e^{-\phi/h} \Delta e^{\phi/h}$ with $\phi(x) = x\cdot \omega$ . Then there exists a Green's function $G_{\Gamma_j}$ defined on $\Omega_j := \Omega \cap \{x_n > f_j(x')\}$ which acts as a bounded operator $G_{\Gamma_j}: L^2(\Omega_j) \to_{h^{-1}} L^2(\Omega_j)$ and $G_{\Gamma_j} : L^{p'}(\Omega_j)\to_{h^{-2}} L^p(\Omega_j)$
%such that $\langle h^2\Delta_\phi^* u, G_\Gamma v \rangle = \langle u, v\rangle$ for all $u\in C^\infty_0(\Omega)$. Furthermore, $G_{\Gamma_j} v\in H^1(\Omega_j)$ for all $v\in L^{p'}(\Omega_j)$  and $G_{\Gamma_j} v \mid_{\Gamma_j} = 0$.
%\end{prop}
To prove Theorem \ref{green's function} in the general case, we first develop the necessary tools for gluing together Green's functions. Let $\Omega$ be a bounded domain and $\Gamma$ be a subset of $\partial\Omega$ which coincides with the graph $\{x_n = f(x')\}$ of a smooth compactly supported function $f$. Without loss of generality we may assume that there is an open neighbourhood $\Omega_{\Gamma} \subset \R^n$ of $\Gamma$ for which $\Omega_{\Gamma} \cap \Omega$ lies in the set $\{x_n > f(x')\}$, and that 
\[
\Omega_{\Gamma} \cap \partial \Omega \cap \{x_n = f(x')\} = \bar\Gamma.
\]
Then $\Gamma' := \Omega_{\Gamma} \cap \partial \Omega$ is an open subset of the boundary such that $\Gamma\subset\subset \Gamma' $ and compact subsets of $\Gamma' \backslash \bar\Gamma$ lies strictly above the graph $x_n = f(x')$.

Let $\chi \in C^\infty_0(\R^n)$ be supported inside $\Omega_\Gamma$ with $\chi = 1$ near $\Gamma$. Then we can arrange that ${\supp} (\chi)\cap \partial\Omega \subset \Gamma'$, and for the derivatives of $\chi$ to have the following support property.
\begin{eqnarray}
\label{support of Dchi}
\exists \epsilon > 0 \mid {\supp}({\bf 1}_{\Omega} D\chi) \subset\{ (x',x_n)\mid x_n \geq f(x') + \epsilon\}.
\end{eqnarray}

In this setting choose an open subset ${\mathcal O}\subset \Omega \cap \{ (x',x_n)\mid x_n> f(x')\}$ which contains $\Gamma'$ as a part of its boundary and whose closure contains the support of $\chi{\bf 1}_\Omega$. Set $G_\Gamma$ to be the Green's function constructed in Proposition \ref{graph green} for the domain ${\mathcal O}$ with vanishing trace on $\Gamma$. We may then define
\begin{eqnarray}\label{Pi estimates}
\Pi_\Gamma : L^{p'}(\Omega) \to_{h^{-2}} L^p(\Omega),\ \ \Pi_\Gamma : L^{2}(\Omega) \to_{h^{-1}} H^1(\Omega)\end{eqnarray}
by 
\[\Pi_\Gamma := \chi {\bf 1}_\Omega (G_\phi - G_\Gamma){\bf 1}_{\mathcal O}.\]
%\lt{do we need this statement about extension? The cutoff $\chi$ already took care of it}
Note that $G_{\Gamma}$ is not defined on the portion of $\Omega$ that lies below the graph of $f$, but this point is rendered moot by the multiplication by $\chi$. Observe that by Proposition \ref{graph green} %$\Pi_\Gamma$, 
one has the trace identity
\begin{eqnarray}
\label{same trace}\Pi_\Gamma v \in H^1(\Omega),\ \ (\Pi_\Gamma v)\mid_{\Gamma} = (G_\phi v)\mid_{\Gamma}, \ \ \ \forall v\in L^{p'}(\Omega).
\end{eqnarray}
\begin{lemma}
\label{remainder of difference in green}
One has the estimates 
\[
h^2\Delta_\phi {\bf 1}_\Omega\Pi_\Gamma {\bf 1}_\Omega : L^{p'}(\Omega)  \to_{h^0} L^2(\Omega),\ \ h^2\Delta_\phi {\bf 1}_\Omega\Pi_\Gamma {\bf 1}_\Omega : L^{2}(\Omega) \to_{h^1} L^2(\Omega).
\]
\end{lemma}
With this lemma we are in a position to construct a general Green's function for the $h^2 \Delta_\phi$ on a general domain $\Omega$. Let $\omega\in \R^n$ be a unit vector and $\Gamma\subset \partial\Omega$ be compactly contained in $\{x\in \partial\Omega \mid \omega \cdot \nu(x)>0\}$ and write $\Gamma$ as a union of its connected components $\Gamma_j$. Without loss of generality we may assume as before that $\omega = (0',1)$. For each $\Gamma_j$ construct $\chi_j$ and $\Pi_{\Gamma_j}$ %{\lt {just noticed that $f$ is used as the function for our graph but simultaneously used as test functions}}
as earlier. One then, by \eqref{same trace}, has that 
\[\Big( G_\phi v- \sum_{j = 1}^k \Pi_{\Gamma_j} v\Big)\mid_\Gamma = 0,\ \ \ \forall v\in L^{p'}(\Omega).\]

Furthermore by Lemma \ref{remainder of difference in green}, $ h^2\Delta_\phi {\bf 1}_\Omega\big(G_\phi - \sum_{j = 1}^k \Pi_{\Gamma_j} \big){\bf 1}_{\Omega}=I+ R'$ with \[R': L^2(\Omega)\to_{h} L^2(\Omega),\ \  R': L^{p'}(\Omega)\to_{h^0} L^{2}(\Omega).\]
Note that as before we can as before invert by Neumann series since $L^{p'}$ gets mapped by $R'$ to $L^2$ with no loss and the Neumann series converge in $L^2$. Theorem \ref{green's function} is now complete by the estimates of \eqref{Pi estimates}, Lemma \ref{sylvester uhlmann}, and Lemma \ref{krs estimate}. %\begin{theorem}
%\label{green with boundary}
%Let $\Omega$ be a bounded set in $\R^n$ with smooth boundary and $\omega\in \R^n$ be a unit vector. If $\Gamma\subset\partial \Omega$ is compactly contained in the set $\{ x\in \partial\Omega \mid \omega \cdot\nu(x)>0\}$. Then there exists a Green's function $G_\Gamma: L^2(\Omega)\to_{h^{-1}} H^1(\Omega)$, and $G_\Gamma : L^{p'}(\Omega)\to_{h^{-2}} L^p(\Omega)$
%such that $\langle h^2\Delta_\phi^* u, G_\Gamma f \rangle = \langle u, f\rangle$ for all $u\in C^\infty_0(\Omega)$. Furthermore, $G_\Gamma f\in H^1(\Omega)$ for all $f\in L^{p'}$  and $G_\Gamma f \mid_{\Gamma} = 0$.
%\end{theorem}
All that remains is to give a proof of Lemma \ref{remainder of difference in green}.

\begin{proof}[Proof of Lemma \ref{remainder of difference in green}]
By Proposition \ref{graph green}, $G_{\Gamma}$ is by construction a right inverse for $h^2\Delta_\phi$ in $\Omega$, and $\chi {\bf 1}_\Omega$ is supported only on $\Omega$, so $\chi h^2 \Delta_\phi {\bf 1}_{\Omega} G_{\Gamma}v(x) = \chi v(x)$ as distributions on $\Omega$. Meanwhile $G_\phi$ is an honest right inverse for $h^2\Delta_\phi$ on $\R^n$, so $h^2\Delta_\phi {\bf 1}_{\Omega} G_\phi = I$ as distributions on $\Omega$. Therefore as distributions on $\Omega$, the only term in $h^2\Delta_\phi \Pi_\Gamma v(x)$ is $[h^2\Delta_\phi, \chi_j] {\bf 1}_\Omega(G_\phi - G_{\Gamma_j}){\bf 1}_{\mathcal O}v(x)$. To analyze this term we will change coordinates by $(x',x_n)\mapsto (x', x_n - f(x'))$ and mark the pushed forward domains, functions and operators with a tilde. Then by the push-forward expression for the operator $G_\Gamma$ stated in Proposition \ref{graph green}, the operator in our term becomes
\[
[h^2\t\Delta_\phi, \t\chi]{\bf 1}_{\t \Omega} (\t G_\phi - (P_s + P_l) {\bf 1}_{\t \Omega} (I + R)){\bf 1}_{\t{\mathcal O} }
\]
where
\[R : L^{p'}(\t\Omega) \to_{h^0} L^2(\t\Omega),\ \ R: L^2(\t\Omega)\to_{h} L^2(\t\Omega).\]
%\[
%S = (1 + hR_{\ell}' + hR_s)
%\]
Computing the commutator $[h^2\t\Delta_\phi, \t\chi]$ explicitly in conjunction with the operator estimates in Proposition \ref{small freq parametrix} and Proposition \ref{estimates for P_l} we have that
%and $(1 + R_l S)$ is inverted by Neumann series. We claim that 
{\Small\begin{equation}\label{LastMultiGraphTerm}
[h^2\t\Delta_\phi, \t\chi]{\bf 1}_{\t \Omega} (\t G_\phi - (P_s + P_l) {\bf 1}_{\t \Omega} S(1 + R_l S)^{-1} ) {\bf 1}_{\t{\mathcal O} } = [h^2\t\Delta_\phi, \t\chi]{\bf 1}_{\t \Omega} (\t G_\phi - (P_s + P_l)) {\bf 1}_{\t \Omega} + E 
\end{equation}}
where
\[
E : L^{p'}(\t \Omega) \to_{h^0} L^2(\t \Omega),\ \ E : L^{2}(\t \Omega) \to_{h^1} L^2(\t \Omega).
\]
%To check this, write
%\[
%S(1 + R_lS)^{-1} = (1 + hR_{\ell}' + hR_s)(1 - R_lS + (R_lS)^2 - \ldots)
%\]
%and recall from the proof of Proposition \ref{graph green} that 
%\[
%R_l : L^{p'}(\Omega) \to_{h^0} L^2(\Omega),\ \ R_l : L^{2}(\Omega) \to_{h^1} L^2(\Omega),
%\]
%while $R_s, R_l' :L^r \rightarrow_{h^0} L^r$ for each $1 < r < \infty$. 
Returning to \eqref{LastMultiGraphTerm}, we see that $E$ has the correct boundedness properties, so it remains only to analyze the first term
\[
[h^2\t\Delta_\phi, \t\chi]{\bf 1}_{\t \Omega} (\t G_\phi - (P_s + P_l)) {\bf 1}_{\t \Omega}.
\]
Since we are only doing the computation in $\t \Omega$, the first order differential operator  $[h^2\t\Delta_\phi, \t\chi]$ commutes with the indicator function ${\bf 1}_{\t \Omega}$, and we have
\[[h^2\t\Delta_\phi, \t\chi]{\bf 1}_{\t \Omega} (\t G_\phi - (P_s + P_l)) {\bf 1}_{\t \Omega}=
{\bf 1}_{\t \Omega} [h^2\t\Delta_\phi, \t\chi](\t G_\phi - (P_s + P_l)) {\bf 1}_{\t \Omega}.
\]
Now $P_s$ maps $L^2$ to $L^2$ with no loss of $h$'s, and $L^{p'}$ to $W^{2,p'} \hookrightarrow_{h^{-1}} H^1$. Meanwhile the commutator $[h^2\t\Delta_\phi, \t\chi]$ maps $H^1$ to $L^2$ with the gain of $h$, so the term involving $P_s$ has the desired behaviour.  Therefore the only term of difficulty is 
\[ [h^2\t\Delta_\phi, \t\chi]{\bf 1}_{\t \Omega} (\t G_\phi - P_l) {\bf 1}_{\t \Omega} = [h^2\t\Delta_\phi, \t\chi]{\bf 1}_{\t \Omega} (I - J^+J)\t G_\phi {\bf 1}_{\t \Omega}  ={\bf 1}_{\t \Omega}[h^2\t\Delta_\phi, \t\chi] J^{-1}{\bf1}_{\R^n_-} J\t G_\phi {\bf 1}_{\t \Omega} \]
By \eqref{support of Dchi} the term ${\bf 1}_{\t \Omega}[h^2\t\Delta_\phi, \t\chi] $ is a first order differential operator whose coefficients are supported in $\{x_n \geq \epsilon >0\}$. The proof then follows from Lemma \ref{disjoint support for J inverse}.
\end{proof}
\end{subsection}

\begin{subsection}{Carleman Estimates} \label{Carleman estimates}

The Carleman estimates in Theorem \ref{carleman estimates} now follow from the existence of the Green's function $G_{\Gamma}$. 

\begin{proof}[Proof of Theorem \ref{carleman estimates}]
Let $u \in C^2(\bar\Omega)$ be a function which vanishes along $\partial\Omega$ and $\partial_\nu u \mid_{\Gamma^c} = 0$, and let $v \in C^{\infty}_0(\Omega)$. Integrating by parts, we have 
\begin{equation}\label{CarlemanIbyP}
\langle h^2\Delta^{*}_{\phi}u, G_{\Gamma}v \rangle_{\Omega} = \langle u, v \rangle_{\Omega} 
\end{equation}
with the boundary terms vanishing because of the boundary conditions on $u$ and the boundary behaviour of $G_{\Gamma}v$. Equation \eqref{CarlemanIbyP} implies that 
\[
\|h^2\Delta_{\phi}u\|_{H^{-1}_{\Gamma}(\Omega)}\|G_{\Gamma}v\|_{H^1_{\Gamma}(\Omega)} \geq |\langle u,v\rangle_{\Omega}|
\]
and 
\[
\|h^2\Delta_{\phi}u\|_{L^{p'}(\Omega)}\|G_{\Gamma}v\|_{L^p(\Omega)} \geq |\langle u,v\rangle_{\Omega}|.
\]
Applying the boundedness results for $G_{\Gamma}$ and taking the supremum over $v \in C^{\infty}_0(\Omega)$ completes the proof. \end{proof}% and the second Carleman estimate will follow if we show that 

\end{subsection}
\end{section}

\begin{section}{Complex Geometrical Optics and the Inverse Problem}
Let $\Omega\subset \R^n$, $\omega \in {\mathbf S}^{n-1}$ and $\Gamma \subset \partial \Omega$ be an open subset of the boundary compactly contained in $\{x\in \partial \Omega\mid \nu(x)\cdot\omega >0\}$ where $\nu_n $ denotes the normal vector. By Theorem \ref{green's function} there exists a Green's function $G_\Gamma$ for $h^2\Delta_\phi$ with vanishing trace on $\Gamma$ and \[G_\Gamma: L^2(\Omega) \to_{h^{-1}} L^2(\Omega),\ \ \ G_\Gamma : L^{p'}(\Omega)\to_{h^{-2}} L^{p}(\Omega).\]

\begin{subsection}{Semiclassical solvability}

Let $\omega$ be a unit vector and $\Gamma\subset \partial\Omega$ be an open subset which is compactly contained in $\{x\in \partial \Omega \mid \nu(x)\cdot\omega >0\}$ we have the following solvability result, resembling the one in \cite{LavNac} (see the explanation of this method in \cite{DosKenSal}), but with an additional term.  
\begin{prop}
\label{solve for rhs}
Let $L \in L^2(\Omega)$ with $\|L\|_{L^2} \leq Ch^2$, and let $q\in L^{n/2}(\Omega)$. For all $a = a_h\in L^\infty$ with $\|a_h\|_{L^\infty} \leq C$, there exists a solution of 
\begin{eqnarray}
\label{solve}
 h^2(\Delta_\phi + q) r = h^2 q a + L\ \ \ r\mid_{\Gamma} = 0\end{eqnarray}
with estimates $\|r\|_{L^2} \leq o(1)$ and $\|r\|_{L^p} \leq O(1)$.
\end{prop}
\proof 
We try solutions of the form $r = G_\Gamma(\sqrt{|q|} v + L)$ for $v \in L^2$ with $\|v\|_{L^2} \leq Ch^2$. Supposing this can be accomplished, then using $\|L\|_{L^2} \leq Ch^2$, \begin{eqnarray*}
\|r\|_{L^2} &\leq& \|G_\Gamma(\sqrt{|q|} v) \|_{L^2 } + \|G_\Gamma(L)\|_{L^2}\\&\leq& \|G_\Gamma(\sqrt{|q|}^\flat v) \|_{L^2 }+ \|G_\Gamma(\sqrt{|q|}^\sharp v) \|_{L^p } + \|G_\Gamma(L)\|_{L^2} \\ &\leq& \frac{C_\epsilon}{h} \|v\|_{L^2} + \frac{C}{h^2} \| \sqrt{|q|}^\sharp v\|_{L^{p'}} + Ch
\end{eqnarray*}
where for any $\epsilon >0$ we decompose $\sqrt{|q|}= \sqrt{|q|}^\sharp + \sqrt{|q|}^\flat$ with $\sqrt{|q|}^\flat \in L^\infty$ and $\|\sqrt{|q|}^\sharp\|_{L^n} \leq \epsilon$.
Therefore, 
\[\|r\|_{L^2} \leq \left(\frac{C_\epsilon} {h} +\frac{C\epsilon}{h^2}\right)\| v\|_{L^{2}} + Ch  = o(1)\]
by taking $h\to 0$ and using that $\|v\|_{^2} \leq Ch^2$. 

For the $L^p$ norm, observe that \[\|L\|_{L^{p'}} \leq \|L\|_{L^2} \leq Ch^2\ \, \, \mbox{ and } \, \, \|\sqrt{|q|} v\|_{p'} \leq \|q\|_{L^{n/2}} \|v\|_{L^2} \leq C h^2.\] The mapping property of $G_\Gamma$ from $L^{p'} \to_{h^{-2}} L^p$ then gives the result.

We now show that we can indeed construct such a $v$. Inserting the ansatz into \eqref{solve} and writing $q = e^{i\theta}|q|$ for some $\theta(\cdot) : \Omega \to \R$ we see that it suffices to construct $v\in L^2$ solving the integral equation
\[(1+ h^2 e^{i\theta} \sqrt{|q|} G_\Gamma\sqrt{|q|}) v = h^2(e^{i\theta} \sqrt{|q|} a - e^{i\theta} \sqrt{|q|} G_\Gamma(L))\]
with $\|v\|_{L^2} \leq Ch^2$. Observe that the right side is $O(h^2)$ in $L^2$ norm due to the fact that $\|L\|_{L^2} \leq Ch^2$ so it suffices to show that $ h^2 e^{i\theta} \sqrt{|q|} G_\Gamma\sqrt{|q|} : L^2\to L^2$ is bounded by $o(1)$ as $h\to 0$ and invert by Neumann series. Indeed, writing $\sqrt{|q|} = \sqrt{|q|}^\sharp + \sqrt{|q|}^\flat$ we have
\[ \sqrt{|q|} G_\Gamma\sqrt{|q|} =  \sqrt{|q|}^\flat G_\Gamma\sqrt{|q|}^\flat+  \sqrt{|q|}^\sharp G_\Gamma\sqrt{|q|}^\flat+ \sqrt{|q|}^{\flat} G_\Gamma \sqrt{|q|}^\sharp.\]
Each of the three pieces have the following mapping properties:
\[\sqrt{|q|}^\flat G_\Gamma\sqrt{|q|}^\flat: L^2 \xrightarrow{\sqrt{|q|}^\flat} L^2 \xrightarrow[h^{-1}] {G_\Gamma} L^2 \xrightarrow{\sqrt{|q|}^\flat} L^2\]
\[\sqrt{|q|}^\sharp G_\Gamma\sqrt{|q|}^\flat: L^2 \xrightarrow{\sqrt{|q|}^\flat} L^2 \hookrightarrow L^{p'} \xrightarrow[h^{-2}] {G_\Gamma} L^p\xrightarrow[o(1)]{\sqrt{|q|}^\sharp} L^2\]
\[\sqrt{|q|}^\flat G_\Gamma\sqrt{|q|}^\sharp: L^2 \xrightarrow[o(1)]{\sqrt{|q|}^\sharp} L^{p'} \xrightarrow[h^{-2}] {G_\Gamma} L^p \xrightarrow{\sqrt{|q|}^\flat} L^p \hookrightarrow L^2\]
Therefore we have that $h^2 e^{i\theta} \sqrt{|q|} G_\Gamma\sqrt{|q|} : L^2\to_{o(1)} L^2$ as $h\to 0$.\qed
\end{subsection}
\begin{subsection}{Ansatz for the Schr\"odinger equation}
We briefly summarize the ansatz construction procedure given in \cite{ksu}; see also the explanation in \cite{Chuthesis}. Let $\phi(x)$ and $\psi(x)$ be linear functions satisfying $D(\phi + i\psi) \cdot D(\phi + i\psi) = 0$. If $\Gamma\subset \partial\Omega$ is an open subset of the boundary satisfying $ D \phi \cdot \nu(x)\geq \epsilon_0 >0$ for all $x\in \bar\Gamma$, we first look to construct a solution to 
\[h^2 \Delta_\phi (e^{i\psi/h}+ a_h) = L,\ \ \ (e^{i\psi/h} + a_h)\mid_{\Gamma} =0\]
with $\|L\|_{L^2} \leq Ch^2$ and $a_h\in L^\infty$. By the fact that $ \nabla \phi \cdot \nu(x) \geq \epsilon_0>0$ for all $x\in \Gamma$, we can apply Borel's lemma to construct $\ell \in C^\infty$ such that
 \[D\ell\cdot D\ell(x) = d(x,\Gamma)^\infty\ \ \ \ell\mid_\Gamma = (\phi+ i\psi)\mid_\Gamma\ \ \ \partial_\nu\ell\mid_\Gamma = -\partial_\nu(\phi+ i\psi)\mid_\Gamma.\]
%Furthermore we can solve the transport equation to infinite order and construct $b\in C^\infty(\Omega)$ such that
%\[h^2 \Delta e^{\ell/h} b = e^{\ell/h}(d(x,\Gamma)^\infty + O_{L^\infty} (h^2))\ \ b\mid_\Gamma = -1\]

Since we are working with linear weights we will need a slightly more general $h$-dependent phase function than $\phi+ i\psi$. Let $\xi \in \R^n$ be a fixed vector which is orthogonal to both $D\phi$ and $D\psi$, and $\psi_h(x)$ be a linear function defined by $\psi_h(x) = (\xi - \omega_h) \cdot x$ where \begin{eqnarray}
\label{omega_h}
\omega_h = \frac{1-\sqrt{1-h^2|\xi|^2}}{h} D\psi 
\end{eqnarray}
 is a vector of length $O(h)$. Observe that in this setting the linear function $\phi + i(\psi +h\psi_h)$ still solves the eikonal equation
\[D(\phi + i(\psi +h\psi_h))\cdot D(\phi + i(\psi +h\psi_h))= 0.\]
We now construct $b\in \C^\infty(\Omega)$ supported close to $\Gamma$ such that
\begin{eqnarray}
\label{solving for amplitude}
e^{-\ell/h} h^2 \Delta (e^{\ell/h} e^{i\psi_h}b) = d(x,\Gamma)^\infty + O_{L^\infty}(h^2) ,\ \ b\mid_\Gamma = -1
\end{eqnarray}
Using the fact that $D\ell \cdot D\ell = d(x,\Gamma)^\infty$ and $D\psi_h = \xi -\omega_h$ with $|\omega_h| \leq Ch$ we see that this amounts to solving the transport equation
\[bD\ell \cdot \xi  + b\Delta\ell + 2D\ell \cdot Db = d(x,\Gamma)^\infty, \ \ \ b\mid_{\Gamma} = -1.\]
Taking advantage of the fact that $-\partial_\nu Re(\ell)\mid_\Gamma = \partial_\nu \phi\mid_\Gamma \geq \epsilon_0>0$ we can again solve the iterative equation and use Borel's Lemma to construct $b\in C^\infty(\Omega)$ supported in an arbitrarily small neighbourhood of $\Gamma$ satisfying this approximate equation. We have therefore constructed $b\in C^\infty$ solving \eqref{solving for amplitude}. 

By the fact that $ \nabla \phi \cdot \nu(x) \geq \epsilon_0>0$ we have, by choosing the support of $b$ sufficiently small, that $Re(\phi(x) - \ell(x))\sim d(x,\Gamma)$ on ${\supp} (b)$. By analyzing separately the case when $d(x,\Gamma) \leq \sqrt{h}$ and $d(x,\Gamma)\geq \sqrt{h}$ we have that \eqref{solving for amplitude} becomes 
\[ h^2 \Delta_\phi (e^{\frac{\ell - \phi}{h}} e^{i\psi_h}b) = O_{L^\infty}(h^2) ,\ \ b\mid_\Gamma = -1\]
By the fact that $h^2 \Delta e^{\frac{\phi+i\psi+hi\psi_h}{h}} = 0$ and $\ell\mid_\Gamma = (\phi+i\psi)\mid_\Gamma$ we have
{\small \begin{eqnarray}
\label{ansatz}
 h^2\Delta_\phi (e^{\frac{i\psi+hi\psi_h}{h}} + e^{\frac{i\psi+hi\psi_h}{h}} a_h) = L,\ \ \|L\|_{L^\infty}\leq Ch^2,\ \ (1 + a_h) \mid_\Gamma = 0.\end{eqnarray}}
where $ a_h := e^{\frac{\ell - \phi - i\psi}{h}} b$ with $\|a_h\|_{L^\infty} \leq C$ and $a_h(x) \to 0$ for all $x\in \Omega$ as $h\to 0$.

%It turns out that the complex phase $\psi$ described here is not general enough for our needs. We need to modify $\psi \mapsto \ \psi_h$ where $\psi_h$ is another linear function satisfying $D(\phi + i\psi_h)\cdot (D\phi+ i\psi_h) = 0$ and differs slightly from $\psi$ in the sense that $\psi_h - \psi = h\xi + \omega  (1- \sqrt{1 - h^2 |\xi|^2})$.
%Observe that we can $D\psi_h \cdot D\psi = D\psi_h \cdot D\phi = 0$.
This discussion allows us to construct the suitable CGO for solving our inverse problem. Indeed, let $\omega$ and $\omega'$ be two unit vectors which are mutually orthogonal. Define $\phi (x) = \omega\cdot x$ and $\psi (x) = \omega' \cdot x$. Let $\xi\in \R^n$ be another vector satisfying $\omega\cdot\xi = \omega'\cdot \xi = 0$ and define $\psi_h(x) := (\xi- \omega_h) \cdot x$ where $\omega_h$ is as in \eqref{omega_h}. Construct $\ell, b\in C^\infty(\Omega)$ so that \eqref{ansatz} is satisfied. Applying Proposition \ref{solve for rhs} to \eqref{ansatz} proves the following
\begin{prop}
\label{CGO}
Let $\omega$ and $\omega'$ be two unit vectors which are mutually orthogonal. Let $\Gamma\subset\partial\Omega$ be an open subset compactly contained in $\{x\in \partial\Omega \mid\omega\cdot \nu(x) >0\}$. For all $q\in L^{n/2}$ there exists solutions to 
\[(\Delta + q) u = 0,\ \ \ u\in H^1(\Omega),\ \ u\mid_{\Gamma} = 0\]
of the form
\[u = e^{\frac{\omega\cdot x + i \omega'\cdot x + hi\psi_h}{h}}(1 + a_h + r)\]
with $\|a_h\|_{L^\infty} \leq C$, $a_h \to 0$ pointwise in $\Omega$ as $h\to 0$. The remainder $r\in L^p$ satisfies the estimates $\|r\|_{L^2} = o(1)$ and $\|r\|_{p} \leq C$ as $h\to 0$.
\end{prop}
%We will also solve the transport equation 
%\[ 2D\ell \cdot D b = -b \Delta \ell,\ \  b\mid_{\Gamma} = -1\]
%to infinite order. {\bf you might need to have this $\Delta \ell$ term on the right hand side at the middle of page 47 of your thesis (instead of zero)}

%The equations for the iterations are,
%\[\sum_{j+k = m} D_t b_j \cdot D_t c_k + \sum_{j+k = m+2} jk b_j c_k = \sum_{j+k = m} c'_k b_j\]
%for some known $c_j$ and $c_j'$ with $c_1\neq 0$ due to the fact that $\partial_\nu \ell = -\partial_\nu \phi >0$ on $\Gamma$. Began by choosing $b_0 = -1$, we see that at each iteration $m$ the only unknown term is $c_1 b_{m+1}$. Using Borel Lemma again we can have $b\in C^\infty$ supported in an arbitrarily small open neighbourhood near $\Gamma$ such that
%\[ 2D\ell \cdot D b = d(x,\Gamma)^\infty -b \Delta \ell,\ \ b\mid_{\Gamma} = -1.\]

%Again, by the fact that $ \nabla \phi \cdot \nu(x) \geq \epsilon_0>0$ we have, by choosing the support of $b$ sufficiently small, that $Re(\phi(x) - \ell(x))\sim d(x,\Gamma)$ on ${\supp} (b)$
%$(\Delta + q) u = 0$ of the form $u = e^{\phi/h} (1 + a_0 + r)$ with $u\mid_\Gamma =0$. We want the remainder $\|r\|_{L^2} = o(1)$ and $\|r\|_{L^p} = O(1)$ as $h\to 0$.

%Without loss of generality assume that $\phi (x) = x_n$. We will construct solution
\end{subsection}
\begin{subsection}{Recovering the Coefficients}
In this section we prove Theorem \ref{main theorem}. Let $\omega$ be a unit vector sufficiently close to $\omega_0$ such that there exists an open set $\Gamma_+$ such that \[\partial\Omega\backslash{\bf B} \subset\subset \Gamma_+\subset\subset \{x\in \partial\Omega \mid \omega \cdot \nu(x) > 0\},\ \ \partial\Omega\backslash{\bf F} \subset\subset \Gamma_-\subset\subset \{x\in \partial\Omega \mid \omega \cdot \nu(x) < 0\}\]
Let $\xi\in \R^n$ be any vector orthogonal to $\omega$ and choose a third vector $\omega'$ of unit length which is perpendicular to both $\xi$ and $\omega$. 

By Theorem \ref{CGO} there exists solutions $u_\pm \in H^1(\Omega)$ solving \[(\Delta + q_1) u_+ = 0,\ \ \ u_+\mid_{\Gamma_+} = 0,\ \ (\Delta+q_2) u_- = 0,\ \ \ u_- \mid_{\Gamma_-} = 0\] of the form
\[u_\pm = e^{\frac{\pm \omega +i\omega' + hi\psi_h^\pm}{h}} (1 + a_h^\pm +r_\pm),\ \ \ \|r_\pm\|_{L^2}= o(1),\ \ \ \|r_\pm\|_{L^p} = O(1)\]
where $\psi_h^\pm(x) := (\pm \xi - \omega') \cdot x$.

Since $u_\pm$ are solutions belonging to $H^1(\Omega)$ and vanish on $\partial\Omega\backslash{\bf B}$ and $\partial\Omega\backslash{\bf F}$ respectively, we have the following boundary integral identity (see Lemma A.1 of \cite{DosKenSal})
\[\int_\Omega \bar u_- (q_1-q_2) u_+ = 0.\]
Inserting the expressions for $u_\pm$ gives
\[0=\int_\Omega e^{2i\xi\cdot x}q(1 + a^-_h a^+_h + a^-_h + a^+_h + a^-_h r_+ + a^+_h r_- + r_- + r_+ + r_+r_-)\]
where $q= q_1 -q_2$. The function $q\in L^{n/2}\subset L^1$ and 
\[\|a_h^\pm\|_{L^\infty}\leq C,\ \ \lim_{h\to 0}a^\pm_h(x) = 0\ \ \forall x\in \Omega\] by \eqref{ansatz}. Therefore, terms $\lim\limits_{h\to0}\int_\Omega e^{2i\xi\cdot x}q( a^-_h a^+_h + a^-_h + a^+_h)  = 0$. For the terms involving $\int_\Omega e^{2i\xi}q a_h^\pm r_\mp$, we note that for all $\epsilon >0$ we may split $q = q^\sharp + q^\flat$ where $q^\flat\in L^\infty$ while $\|q^\sharp\|_{L^{n/2}} \leq \epsilon$. Then, using the fact that $\|a_h^\pm\|_{L^\infty}\leq C$,
\[\left|\int_\Omega e^{2i\xi\cdot x}q a_h^\pm r_\mp \right| \leq C( \|q^\flat\|_{L^\infty} \| r_\mp\|_{L^2} + \|q^\sharp\|_{L^{n/2}} \|r_\mp\|_{L^p})\]
where $p = \frac{2n}{n-2}$. By the estimates on $r_\mp$ given in Proposition \ref{CGO} we have that $\lim\limits_{h\to 0}\|r_\mp\|_{L^2} = 0$ and $\|r_\mp\|_{L^p} \leq C$. Therefore, the limit
\[\lim_{h\to 0}\left|\int_\Omega e^{2i\xi\cdot x}q a_h^\pm r_\mp \right| \leq C\epsilon\]
for all $\epsilon >0$ and therefore the limit vanishes. The terms $\int_\Omega e^{2i\xi}q(r_- + r_+)$ can be estimated the same way.
For the last term, we again decompose, for all $\epsilon >0$, $q= q^\flat + q^\sharp$. The integral $|\int_\Omega e^{2i\xi\cdot x} qr_-r_+|$ is then estimated by
\[\int_\Omega| q^\flat r_-r_+| + \int |q^\sharp r_-r_+| \leq \|q^\flat\|_{L^\infty}\|r_+\|_{L^2} \|r_-\|_{L^2} + \|q^\sharp\|_{L^{n/2}} \|r_-\|_{L^p} \|r_+\|_{L^p}\] 
The $L^p$ norms of $r_\pm$ stay uniformly bounded while the $L^2$ norms vanish when $h\to 0$. Therefore the limit
\[\lim_{h\to 0}\left| \int_\Omega e^{2i\xi\cdot x} qr_-r_+ \right| \leq C\|q^\sharp\|_{L^{n/2}} \leq C\epsilon\]
for all $\epsilon >0$ and therefore vanishes.

This means that ${\mathcal F}(q) (\xi) = 0$ for all $\xi$ which are orthogonal to $\omega$. Note that varying $\omega$ in a small neighbourhood does not change the fact that $\Gamma$ lies in the set $\{x \in \partial \Omega|\omega \cdot \nu(x) > 0\}$, and so the construction in Proposition \ref{CGO} still applies. Then varying $\omega$ in a small neighbourhood and using the analyticity of the Fourier transform for $q$ compactly supported we have that $q = q_1 - q_2= 0$. \qed 
\end{subsection}
\end{section}

\section{Appendix}

Here we will provide proofs for Proposition \ref{sobolev mapping} and Proposition \ref{sobolev composition} from Section 2.  

To begin, suppose $a\in S^{k}_0(\R^n)$ be a symbol whose spatial dependence is in $x'$ only and compactly supported. We then have the following expression for the quantization of their product:
{\small \begin{eqnarray}
\label{mixed composition expression}
a(x',hD) f =  \int e^{i\lambda'\cdot x'} \int e^{-i\lambda'\cdot z'} \int e^{i\xi\cdot x} \frac{ (1 +\Delta_{z'})^Na(z', h\xi)}{(1 + |\lambda'|^2)^N} {\mathcal F}(u)(\xi) d\xi dz' d\lambda'
\end{eqnarray}}

\begin{prop}
\label{Lr composition of mixed symbols}
Let $a(x', \xi)$ be in $S^{0}_1(\R^n)$ or $S^{-k(n)}_0(\R^n)$ for some $k(n)$ large depending only on the dimension. Suppose  $a(x', \xi)$ depends only on $x'\in \R^{n-1}$ in the spatial variable and $D_{x'}a(x',\xi) = 0$ if $x'$ is outside of a fixed compact set.  If $b(x',\xi')\in S^{0}_1(\R^{n-1})$ with $D_{x'} b(x', \xi') = 0$ outside of a compactly supported set, then \[ab(x', hD) : L^r \to L^r\]
with norm {\tiny \[\|ab(x', hD)\|_{L^r\to L^r} \leq C \sup\limits_{\substack{z', \xi,\\ |\alpha| \leq k(n)}} |(1+ \Delta_{z'})^N \partial_\xi^\alpha a(z', \xi)| \langle \xi\rangle ^{|\alpha|}\sup\limits_{\substack{z', \xi,\\ |\alpha| \leq k(n)}} |(1+ \Delta_{z'})^N \partial_{\xi'}^\alpha b(z', \xi')| \langle \xi'\rangle ^{|\alpha|}\]}
where the constant $C$ depends linearly on the volume of the support of $D_{x'}ba(x',\xi)$ in $x'$ and $N$ depends only on the dimension.
\end{prop}
\proof In the constant coefficient case this is a direct consequence of Mihlin's multiplier theorem applied first to all variable then to $\xi'$ variables. We can therefore assume without loss of generality that either $b(x', \xi) = 0$ or $a(x', \xi) = 0$ for $x'$ outside of a fixed compact set.

Then apply Minkowski to expression \eqref{mixed composition expression} for $N$ chosen to be large enough and for each $z'\in {\supp}(ba(\cdot,\xi))$ apply the constant coefficient estimate for Fourier multipliers on $L^r$. \qed

An immediate Corollary is the mapping property from Sobolev spaces :
\begin{cor}
%\label{sobolev mapping}
If $a(x', \xi)\in S^{k}_1 S^{\ell}_1 \cup S_1^{k} S^{-k(n) +\ell}_0$ then
\[a(x',hD):W^{k,r}W^{\ell,r} \to L^r\] with norm uniformly bounded in $h$.  \end{cor}
\proof
Pre-composition yields that $a(x',hD) \langle hD'\rangle^{-k} \langle hD\rangle^{-\ell}$ is a quantization of a symbol in $S^{0}_1 S^{0}_1 \cup S^{0}_1 S^{-k(n)}_0$ and therefore takes $L^r \to L^r$. This shows that \[a(x',hD) : W^{k,r}W^{\ell,r} \to L^r\ \ \forall 1<r<\infty.\]\qed

Composition of two $\Psi$DO operators in this class can be described by the composition calculus $ b(x, hD)a(x, hD)  = ab(x, hD) + h \sum\limits_{|\alpha| = 1}(\partial_\xi^\alpha b \partial_x^\alpha a)(x, hD) +h^2 m(x,hD)$ with the remainder explicitly computed as
{\tiny \begin{eqnarray}
\label{composition remainder}
\ \ \ \ \ \ \ m(x,\xi) = \sum\limits_{|\alpha|=2} \int_{\R^{2n}} \frac{e^{i\eta\cdot y}}{\langle \eta\rangle^{2N}\langle y\rangle^{2N}} (I +\Delta_\eta)^{N}\partial^\alpha_\eta b(x, \eta + \xi)(I + \Delta_y)^{N} \int_0^1 \partial_x^\alpha a (x+\theta h y, \xi)d\theta dyd\eta\ \ \forall N\in \mathbb N
\end{eqnarray}} 
This leads to the following statement about the remainder term of the composition:
\begin{lemma}
\label{Lr composition formula}
Let $a \in S^{k_1}_1 S^{\ell_1}_1 \cup S^{k_1}_1 S^{-k(n)+\ell_1}_0$ and $b \in S^{-k_1}_1 S^{-\ell_1}_1 \cup S^{-k_1}_1 S^{-\ell_1 -k(n)}_0$ then one has $ b(x', hD)a(x', hD)  = ab(x', hD) + h \sum\limits_{|\alpha| = 1}(\partial_\xi^\alpha b \partial_x^\alpha a)(x', hD) +h^2 m(x',hD)$ with $m(x',hD) : L^r\to L^r$ norm independent of $h>0$. 
\end{lemma}
\proof We have that 
\[ b(x'hD) a(x', hD) = ab(x',hD) + h\sum\limits_{|\alpha|=1} \partial^\alpha_{\xi} b \partial^\alpha_{x'} a (x',hD) + h^2m(x',hD)\]
where $m(x,\xi)$ is given by \eqref{composition remainder}. By taking $N$ large enough in $\eqref{composition remainder}$ we see that 
\begin{eqnarray}
\label{superposition of remainder}
m(x', hD) u = \sum\limits_{|\alpha| = 2} \int_{\R^{2n}} \frac{e^{iy\cdot \eta}}{\langle \eta\rangle^{N}\langle y\rangle^{N}}\int_0^1 m^{\alpha,j}_{\theta,y,h,\eta}(x', hD) u  d\theta dy d\eta
\end{eqnarray}
where for each $(\alpha,\theta,y,h,\eta)$, $m^{\alpha}_{\theta,y,h,\eta} (x', \xi) \in S^0_1 S^0_1 \cup S^0_1 S^{-k(n)}_0$ is a symbol of the form $m^{\alpha}_{\theta,y,h,\eta} (x', \xi) =\langle \eta\rangle^{-N}\langle y\rangle^{-N}(I +\Delta_\eta)^N  \partial^\alpha_\eta b(x', \eta+ \xi)(1+ \Delta_y)^N  \partial_x^\alpha a(x + h\theta y, \xi) $. Since $a\in S^{k_1}_1S_1^{\ell_1} \cup S^{k_1}_1 S^{-k(n) +\ell_1}_0$ and $b\in S^{-k_1}_1S_1^{-\ell_1} \cup S^{-k_1}_1 S^{-k(n) -\ell_1}_0$ we may write $a = a^ta^v$ and $b= b^tb^v$ where \[a^t(x',\xi')\in S^{k_1}_1, b^v(x',\xi') \in S^{-k_1}_1, a^v(x',\xi) \in S_1^{\ell_1}\cup S^{-k(n) +\ell_1}_0, b^v(x',\xi)\in S_1^{-\ell_1}\cup S^{-k(n)-\ell_1}_0.\]
We see then that for each $(\alpha,\theta,y,h,\eta)$ the symbol $m^{\alpha}_{\theta,y,h,\eta} (x', \xi)$ consists of finitely many (depending on the choice of $N$) terms of the from
{\small\[\langle \eta\rangle^{-N}\langle y\rangle^{-N} \partial_\xi^{\beta_1} b^t_\eta (x',\xi) (\theta h)^{|\beta_2|} \partial_{x'}^{\beta_2} a^t_{\theta, h y'}(x',\xi) \partial_\xi^{\beta_3} b^v_\eta (x',\xi) (\theta h)^{|\beta_4|} \partial_{x'}^{\beta_4} a^v_{\theta, h y'}(x',\xi)\]}
which is a symbol in $S^0_1 S^0_1 \cup S^0_1 S^{-k(n)}$. Here $b_\eta (x',\xi) := b(x', \xi+ \eta)$ and $a_{\theta, h y'}(x',\xi) := a(x'+ \theta h y',\xi)$
Applying Proposition \ref{Lr composition of mixed symbols} to each of these terms and choosing $N \geq k(n)$ we have that 
\[\sup\limits_{(\alpha,\theta,y,h,\eta)} \|m^{\alpha}_{\theta,y,h,\eta}(x',hD)\|_{L^r\to L^r} \leq C_N\]
Here we used Peeter's inequality $\frac{\langle \xi\rangle}{\langle \xi +\eta\rangle \langle \eta\rangle} \leq C$.
Choosing $N\geq n+2$ in \eqref{superposition of remainder} we get that 
{\small\[\|m(x',hD)\|_{L^r\to L^r} \leq C\sup\limits_{(\alpha,\theta,y,h,\eta)} \|m^{\alpha}_{\theta,y,h,\eta}(x',hD)\|_{L^r\to L^r} \int_{\R^{2n}}  \langle y\rangle^{-N} \langle \eta\rangle^{-N} d\eta dy . \]}\qed

The composition formula given by Lemma \ref{Lr composition formula}
 in conjunction with the mapping property asserted in Proposition \ref{Lr composition of mixed symbols}
 also allows us to deduce, Proposition \ref{sobolev mapping} by composition with suitable powers of $\langle hD'\rangle \langle hD\rangle$. 
\begin{prop}[Proposition \ref{sobolev mapping}]
If $b(x',\xi')\in S^{k}_1$ and $a(x',\xi) \in S^{\ell}_1 \cup S^{-k(n)+\ell}_0$ then \[ba(x', hD) : W^{m,r} W^{l,r} \to W^{m-k,r} W^{l - \ell,r}\] with norm
{\tiny \[\|ab(x', hD)\| \leq C \sup\limits_{\substack{z', \xi,\\ |\alpha| \leq k(n)}} |(1+ \Delta_{z'})^N \partial_\xi^\alpha a(z', \xi)| \langle \xi\rangle ^{|\alpha|-\ell}\sup\limits_{\substack{z', \xi,\\ |\alpha| \leq k(n)}} |(1+ \Delta_{z'})^N \partial_{\xi'}^\alpha b(z', \xi')| \langle \xi'\rangle ^{|\alpha|-k}\]}

\end{prop}
\proof
Since pre-composition by $\langle hD'\rangle^{-k}\langle hD\rangle^{-\ell}$ amounts to multiplication of symbols without remainders, it suffices to show that symbols $a(x',\xi)\in S^{k}_1 S^{\ell}_1 \cup S^{k}_1 S^{-k(n)+\ell}_0$ take $L^r \to W^{-k,r} W^{-\ell,r}$. Indeed, by Lemma \ref{Lr composition formula} we have that \[\langle hD'\rangle^{-k}\langle hD\rangle^{-\ell} a(x',hD) = c(x',hD) + h^2 m(x',hD)\] where $c(x',hD) = \langle \xi'\rangle^{-k}\langle \xi\rangle^{-\ell} a(x',\xi) + h\sum\limits_{|\alpha|=1} \partial_{\xi}^\alpha( \langle \xi'\rangle^{-k}\langle \xi\rangle^{-\ell}) \partial^\alpha_{x'} a(x',\xi)$ and $m(x',hD): L^r\to L^r$.

To estimate the operator norm by the size of the symbol, using \eqref{composition remainder} and estimate the remainder as in the proof of Proposition \ref{Lr composition of mixed symbols}.\qed

Now we turn to the proof of Proposition \ref{sobolev composition}.
\begin{prop}[Proposition \ref{sobolev composition}]
If $a\in S^{k_1}_1S_1^{\ell_1} \cup S^{k_1}_1 S^{-k(n) +\ell_1}_0$ and $b\in S^{k_2}_1S_1^{\ell_2} \cup S^{k_2}_1 S^{-k(n) +\ell_2}_0$ then
\[  b(x'hD)a(x',hD) = ab(x',hD) + h\sum\limits_{|\alpha|=1} (\partial_{\xi}^\alpha b\partial_{x'}^\alpha a)(x',hD) + h^2 m(x',hD)\]
where $m(x',hD): W^{k,r}W^{\ell,r} \to W^{k - k_1-k_2,r}W^{\ell-\ell_1-\ell_2,r}$. 
\end{prop}
\proof 
The proof is exactly the same as Proposition \ref{Lr composition of mixed symbols} except that to show the boundedness of the remainder in the mixed Sobolev norms one uses Proposition \ref{sobolev mapping}.

We have that 
\[ b(x'hD) a(x', hD) = ab(x',hD) + h\sum\limits_{|\alpha|=1} \partial^\alpha_{\xi} b \partial^\alpha_{x'} a (x',hD) + h^2m(x',hD)\]
where $m(x,\xi)$ is given by \eqref{composition remainder}. By taking $N$ large enough in $\eqref{composition remainder}$ we see that 
\begin{eqnarray}
\label{superposition of remainder 2}
m(x', hD) u = \sum\limits_{|\alpha| = 2} \int_{\R^{2n}} \frac{e^{iy\cdot \eta}}{\langle \eta\rangle^{N}\langle y\rangle^{N}}\int_0^1 m^{\alpha,j}_{\theta,y,h,\eta}(x', hD) u  d\theta dy d\eta
\end{eqnarray}
where for each $(\alpha,\theta,y,h,\eta)$, $m^{\alpha}_{\theta,y,h,\eta} (x', \xi) \in S^0_1 S^0_1 \cup S^0_1 S^{-k(n)}_0$ is a symbol of the form $m^{\alpha}_{\theta,y,h,\eta} (x', \xi) =\langle \eta\rangle^{-N}\langle y\rangle^{-N}(I +\Delta_\eta)^N  \partial^\alpha_\eta b(x', \eta+ \xi)(1+ \Delta_y)^N  \partial_x^\alpha a(x + h\theta y, \xi) $. Since $a\in S^{k_1}_1S_1^{\ell_1} \cup S^{k_1}_1 S^{-k(n) +\ell_1}_0$ and $b\in S^{-k_1}_1S_1^{-\ell_1} \cup S^{-k_1}_1 S^{-k(n) -\ell_1}_0$ we may write $a = a^ta^v$ and $b= b^tb^v$ where \[a^t(x',\xi')\in S^{k_1}_1, b^v(x',\xi') \in S^{-k_1}_1, a^v(x',\xi) \in S_1^{\ell_1}\cup S^{-k(n) +\ell_1}_0, b^v(x',\xi)\in S_1^{-\ell_1}\cup S^{-k(n)-\ell_1}_0.\]
We see then that for each $(\alpha,\theta,y,h,\eta)$ the symbol $m^{\alpha}_{\theta,y,h,\eta} (x', \xi)$ consists of finitely many (depending on the choice of $N$) terms of the from
{\small\[\langle \eta\rangle^{-N}\langle y\rangle^{-N} \partial_\xi^{\beta_1} b^t_\eta (x',\xi) (\theta h)^{|\beta_2|} \partial_{x'}^{\beta_2} a^t_{\theta, h y'}(x',\xi) \partial_\xi^{\beta_3} b^v_\eta (x',\xi) (\theta h)^{|\beta_4|} \partial_{x'}^{\beta_4} a^v_{\theta, h y'}(x',\xi)\]}
which is a symbol in $S^0_1 S^0_1 \cup S^0_1 S^{-k(n)}$. Here $b_\eta (x',\xi) := b(x', \xi+ \eta)$ and $a_{\theta, h y'}(x',\xi) := a(x'+ \theta h y',\xi)$.

Applying Proposition \ref{Lr composition of mixed symbols} to each of these terms and choosing $N \geq k(n)$ we have that 
\[\sup\limits_{(\alpha,\theta,y,h,\eta)} \|m^{\alpha}_{\theta,y,h,\eta}(x',hD)\|_{L^r\to L^r} \leq C_N\]
Here we used Peeter's inequality $\frac{\langle \xi\rangle}{\langle \xi +\eta\rangle \langle \eta\rangle} \leq C$.
Choosing $N\geq n+2$ in \eqref{superposition of remainder} we get that 
{\small\[\|m(x',hD)\|_{L^r\to L^r} \leq C\sup\limits_{(\alpha,\theta,y,h,\eta)} \|m^{\alpha}_{\theta,y,h,\eta}(x',hD)\|_{L^r\to L^r} \int_{\R^{2n}}  \langle y\rangle^{-N} \langle \eta\rangle^{-N} d\eta dy  \]} \qed

\end{document}